\pdfoutput=1
 
\documentclass{article}
\usepackage{amsmath}
\usepackage[utf8,ansinew]{inputenc}
\usepackage{times}
\usepackage{units}
\usepackage{epic}
\usepackage{eepic}
\usepackage{amssymb}
\usepackage{graphicx}
\usepackage{subfigure}
\usepackage{algpseudocode}
\usepackage{algorithm}
\usepackage{graphicx}
\usepackage{caption}
\usepackage{graphicx,subfigure}
\usepackage{xfrac}
\usepackage{moreverb}
\usepackage[multiple]{footmisc}

\newcommand\BibTeX{{\rmfamily B\kern-.05em \textsc{i\kern-.025em b}\kern-.08em
T\kern-.1667em\lower.7ex\hbox{E}\kern-.125emX}}
\let\oldReturn\Return
\renewcommand{\Return}{\State\oldReturn}
\let\oldCall\Call
\renewcommand*\Call[2]{\State \oldCall{#1}{#2}}

\newcommand{\vertiii}[1]{{\left\vert\kern-0.25ex\left\vert\kern-0.25ex\left\vert #1 
    \right\vert\kern-0.25ex\right\vert\kern-0.25ex\right\vert}}
\providecommand{\e}[1]{\ensuremath{\times 10^{#1}}}

\newtheorem{lemma}{Lemma}
\newtheorem{problem}{Problem}
\newtheorem{discretization}{Discretization}

\def\nat{\mathbb{N}}
\def\real{\mathbb{R}}

\def\cA{{\cal A}}
\def\cM{{\cal M}}

\def\cD{{\cal D}}

\def\cB{{\cal B}}
\def\cU{{\cal U}}

\def\cR{{\cal R}}

\def\bt{{\mathbf{t}}}
\def\bi{{\mathbf{i}}}
\def\bI{{\mathbf{I}}}
\def\bx{{\mathbf{x}}}

\def\bXi{{\mathbf{\Xi}}}
\def\slin{{\mbox{\scriptsize lin}}}
\def\sPre{{\mbox{\scriptsize prew}}}
\def\sFull{{\mbox{\scriptsize full}}}
\def\sspan{{\mbox{span} \, }}
\def\semiortho{{ \mbox{\scriptsize semi-ortho} }}
\def\coa{{ \mbox{\scriptsize co} }}
\def\dim{{ \mbox{dim} }}
\def\res{{ \mbox{\scriptsize res} }}
\def\prol{{ \mbox{\scriptsize prol} }}
\def\Id{{ \mbox{Id} }}

\newcommand{\Markiere}[1]{{\bf #1}}

\def\setto{~~~{\raisebox{-0.1cm}{\scriptsize set}  \atop \raisebox{0.1cm}{$\longleftarrow$} } ~~~}

\def\volumeyear{2015}

\begin{document}

\title{A Sparse Grid Discretization with Variable Coefficient in High Dimensions}

\author{Rainer~Hartmann and Christoph~Pflaum}



\maketitle

\begin{abstract}
We present a Ritz-Galerkin discretization on sparse grids using pre-wavelets, which allows to
solve  elliptic differential equations with variable coefficients for dimension $d=2,3$ and higher dimensions
$d>3$. The method applies multilinear finite elements.
We introduce an efficient algorithm for matrix vector multiplication using a Ritz-Galerkin discretization 
and semi-orthogonality.

This algorithm is based on standard 1-dimensional restrictions and prolongations, a simple pre-wavelet stencil, and the classical operator dependent stencil for multilinear finite elements.
Numerical simulation results are presented for a 3-dimensional problem on a curvilinear bounded domain
and for a 6-dimensional problem with variable coefficients.

Simulation results show a convergence of the discretization according to the approximation properties of the
finite element space. The condition number of the stiffness matrix can be bounded below $10$ 
using a standard diagonal preconditioner.
\end{abstract}




\section{Introduction}
\vspace{-2pt}

A finite element discretization of an elliptic symmetric PDE calculates the best approximation
with respect to the energy norm. Since finite element method uses polynomials to construct
a finite element space, the convergence of such a method can be proven using Strang\textquoteright s lemma or Lax-Milgram theorem.
However, difficulties arise in the application to high dimensional problems. Then the computational
amount increases by $O(N^d)$, where $N$ is the number of grid points in one direction and
$d$ is the dimension of the space. This exponential growth of the computational amount restricts
the application of the finite element method to dimensions $d\leq3$. A suggestion to solve this problem is 
to use sparse grids (see \cite{ZengerSparseGrids}).
With sparse grids, one can construct a subspace of the classical finite element spaces on full grids.
The dimension of this subspace reduces to $O(N (\log N)^{d-1})$.
However, solving the resulting linear equation efficiently is a difficult task.
Suitable algorithms were constructed for constant coefficients and cubical domains (see \cite{bungartz92duenne} and \cite{balder96solution}).
Those algorithms evaluate the matrix vector multiplication in an efficient way.
Therefore, using a Ritz-Galerkin discretization,
one matrix vector-multiplication can be performed by $O(N (\log N)^{d-1})$ operations.
However, a variable coefficient in the operator leads to $O(N^d)$ operations. 
So far, sparse grids have had very limited range of applications since most partial differential equations
in natural science or engineering science include variable coefficients.

The first sparse grid discretization with variable coefficients, using a Ritz-Galerkin approach,
is presented in \cite{Pfl11}. This discretization applies the semi-orthogonality property of 
standard hierarchical basis functions (see Section \ref{SecSGdis}).
A complete convergence theory is given in \cite{Pfl8} for a 2-dimensional case.
The discretization leads to a symmetric stiffness matrix in case of a symmetric bilinear form.
Nevertheless, an extension to higher dimension problems is not possible for standard hierarchical basis
functions since hierarchical basis functions do not satisfy a semi-orthogonality property for $d\geq3$.

Other discretizations of PDEs with variable coefficients are presented in \cite{achatz03higher} and
\cite{griebelDiff}. The discretization in \cite{achatz03higher} can be treated as a finite element discretization,
while the discretization in \cite{griebelDiff} is a finite difference discretization. Both discretizations lead 
to a non-symmetric linear equation system for symmetric problems, which is an undesired property. 
Additionally, a convergence proof is missing for both discretizations.
Therefore, convergence of these methods is not guaranteed in higher dimensions. 
Furthermore, the discretization in \cite{achatz03higher} requires high order interpolation operators, which
increases the computational amount. However, simulation results presented
in literature show an optimal convergence for certain 2-dimensional and 3-dimensional problems.

In this paper, we present a new method to discretize partial elliptic differential equations
on sparse grids (see Section \ref{SecSGdis}).
This discretization uses pre-wavelets and their semi-orthogonality property (see \cite{Pfl11}).
It is well-known that pre-wavelets and wavelets can be used to discretize partial differential equations
(see \cite{Osw}, \cite{UrbanDahmen}, \cite{VassiWang}, and \cite{Pfl13}).
In the context of sparse grids, they can even lead to natural discretizations of elliptic partial differential equations with variable coefficients.
The elementary convergence theory of such discretizations is presented for a Helmholtz problem in \cite{PflaumHartmann}.
In Section \ref{SectionMatrixVector}, we present an algorithm that efficiently evaluates
the matrix vector multiplication with the discretization matrix. The algorithm applies only
 standard 1-dimensional restriction and prolongation operators, a simple pre-wavelet stencil of size $5$,
and the classical stencil operator for multilinear finite elements.
This operator dependent stencil is a well-known $9$-point stencil for bilinear elements and
a $27$-point stencil for trilinear elements. However, in the 6-dimensional
case the size of this stencil increases to $729=3^6$.
The difficulty of the algorithm is to apply all operators in the correct sequential ordering.

In Section \ref{simResults} simulation results are presented 
for the 3-dimensional Poisson\textquoteright s problem on a curvilinear bounded domain
and for a 6-dimensional Helmholtz problem with a variable coefficient. The simulation result for 
Poisson\textquoteright s problem implies that sparse grids are not restricted to cubical domains.
To our knowledge, the numerical result for the 6-dimensional Helmholtz equation is the first simulation result
for an 6-dimensional Ritz-Galerkin finite element discretization of a elliptic PDE with variable coefficients.

This paper is restricted to non-adaptive grids. However, it is clear that the algorithm
can be extended to adaptive sparse grids as shown in \cite{Pfl11}.

\section{Sparse Grid Discretization}
\vspace{-2pt}
\label{SecSGdis}

\begin{figure}[hbt]
 \centering
   \includegraphics[width=0.45\textwidth]{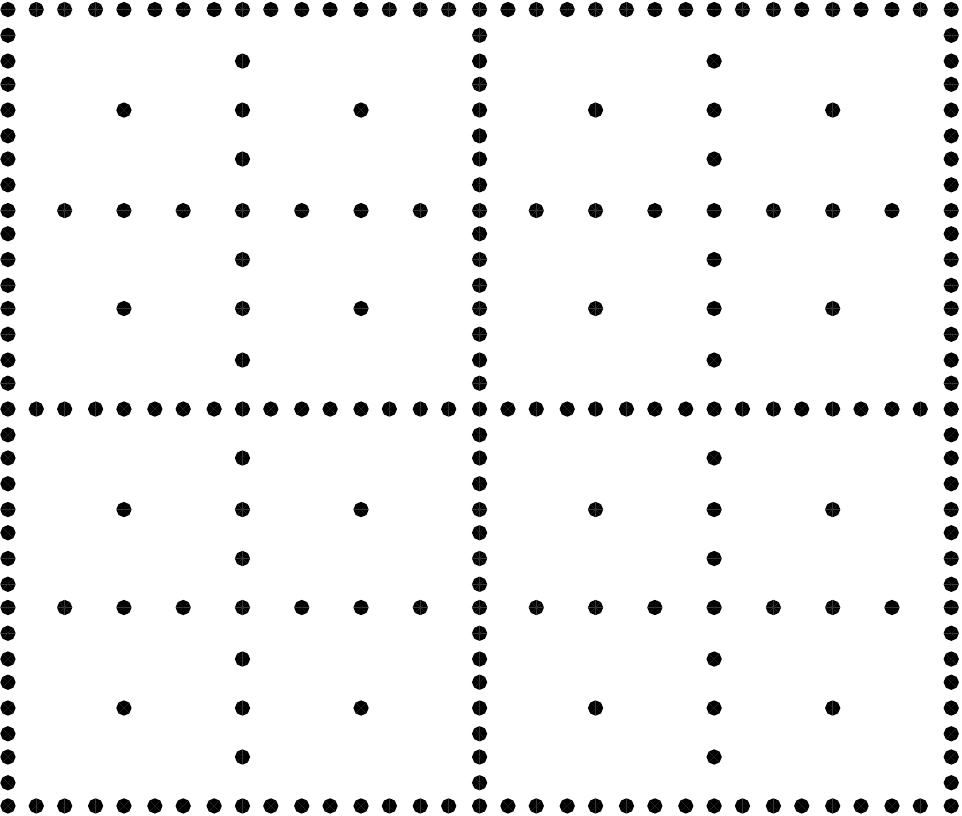}
 \caption{Example of a 2-dimensional sparse grid $\cD_n$.}
    \label{sparsegrid2D}
\end{figure}

Let $d \geq 2$ be the dimension of space and $\Omega = [0,1]^d$. 
Consider the following elliptic differential equation:

\begin{problem}
\label{PropOrg}
 Let $f \in L^2(\Omega)$, $A \in (L^\infty(\Omega))^{d \times d}$ and $\kappa \in L^\infty(\Omega), \kappa \geq 0$ be given.
 Furthermore, assume that $A$ is symmetric and uniformly positive definite. This means that there is a
 $\alpha>0$ such that $v^T A(\bx) v > \alpha v^T v$ for almost every $\bx \in \Omega$ and every vector $v \in \real^d$.
Find $u \in H_0^1(\Omega)$ such that
\[
  \int_\Omega (\nabla u)^T A(\bx) \nabla v + \kappa(\bx) uv ~d\bx = \int_\Omega f v_h ~d\bx \quad \forall v \in H_0^1(\Omega).
\]
\end{problem}
Our aim is to find an efficient sparse grid finite element discretization that can
be used even for large dimension $d$. A typical 2-dimensional sparse grid
is depicted in Figure \ref{sparsegrid2D} and a 3-dimensional sparse grid in Figure \ref{sparsegrid3D}.

\begin{figure}%
  \begin{minipage}{0.45\textwidth}
    \includegraphics[width=\textwidth]{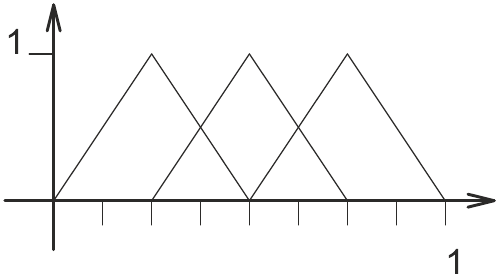}
    \caption{1-dimensional nodal basis functions.}
    \label{nodalLin}
  \end{minipage}%
  \qquad
  \begin{minipage}{0.45\textwidth}
    \includegraphics[width=\textwidth]{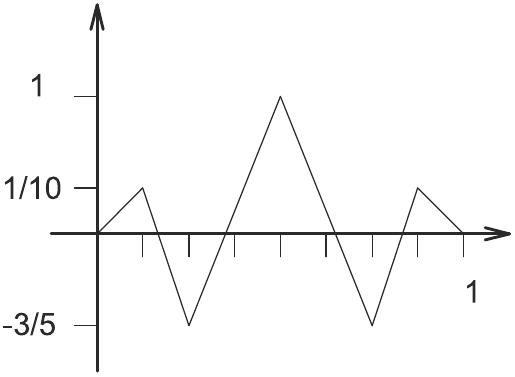}
    \caption{1-dimensional pre-wavelet function.}
    \label{prewavelets}
  \end{minipage}
\end{figure}

Finite elements on sparse grids are constructed by tensor products of 1-dimensional finite elements.
Here, we apply piecewise linear elements in 1D. Let us explain the construction of the  sparse grid finite element
space in more detail. To this end, define the 1-dimensional grid
\begin{eqnarray*}
 \Omega_k & = & \{ 2^{-k-1} i ~|~ i \in I_k \}, \\
 I_k & = & \{ i ~|~ i = 1, ..., 2^{k+1}-1 \},
\end{eqnarray*}
for $k \in \nat_0$, where $I_k$ is the index set of $\Omega_k$. Observe that
$\Omega_0 \subset \Omega_1 \subset \Omega_2 \subset ...$ .
The complementary index set is defined by
\[
 \Xi_t = I_t \backslash I_{t-1}.
\] 
Now, let $V_k$ be the space of piecewise linear functions of mesh size $2^{-k-1}$
and $v^\slin_{k,i}$ the corresponding nodal basis function at point $2^{-k-1} i \in \Omega_k$ 
(see Figure \ref{nodalLin}).

Using these functions, we define pre-wavelets $\varphi_{t,i}$, $i \in \Xi_t$ by (see Figure \ref{prewavelets})
\[
 \varphi_{t,i} =
\left\{
\begin{array}{lcl}
\frac{9}{10} v^\slin_{t,i} - \frac{3}{5} v^\slin_{t,i+1}  + \frac{1}{10} v^\slin_{t,i+2}  & \mbox{if} & 1 = i \\[2ex]
 v^\slin_{t,i} - \frac{3}{5} (v^\slin_{t,i+1} + v^\slin_{t,i-1}) + \frac{1}{10} (v^\slin_{t,i+2} + v^\slin_{t,i-2}) & \mbox{if} & 3 \leq i \leq 2^{t}-3 \\[2ex] 
\frac{9}{10} v^\slin_{k,i} - \frac{3}{5} v^\slin_{t,i-1}  + \frac{1}{10} v^\slin_{t,i-2}  & \mbox{if} & i = 2^{t}-1
\end{array}
\right. ,
\]
for $t \in \nat$ and $\varphi_{0,1} = v^\slin_{0,1}$.
An important property of these functions is the $L^2$-orthogonality for pre-wavelets of different levels
\begin{equation}
\label{orthoPrewave}
\int_0^1 \, \varphi_{t,i} \, \varphi_{t',i'} ~dx \, = \,0 \quad \mbox{if $t\not= t'$}.
\end{equation}

\begin{figure}
 \centering
 \includegraphics[width=0.7\textwidth]{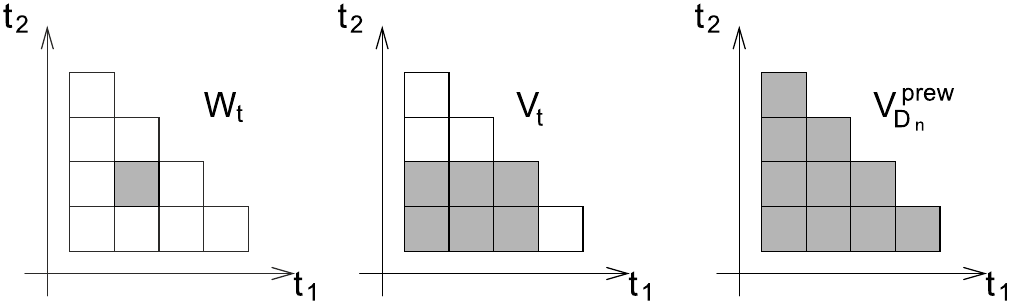}
 \caption{Example of spaces $W_{\bt} = W_{1,1}$, $V_\bt^\sFull = V_{2,1}^\sFull$, and $V_{\cD_n}^\sPre = V_{\cD_3}^\sPre$.}
 \label{Picspaces}
\end{figure}

Let us introduce the following abbreviations for a multi-index $\bt=(t_1,...,t_d) \in \nat_0^d$:

\begin{eqnarray*}
| \bt | & := & \sum_{i=1}^d |t_i|, \\
\bt \leq \bt' & \mbox{~if~} &   t_i \leq t'_i \quad \forall i=1,...,d, \\
\max(\bt,\bt') & := & ( \max(t_1,t_1'), ... , \max(t_d,t_d') ),  \quad \mbox{and}   \\
\max(\bt)    & := & \max(t_1,...,t_d) .
\end{eqnarray*}
Using these abbreviations, the sets of tensor product indices are defined
\begin{eqnarray*}
   \bI_\bt & = & \Big\{ (i_1,...,i_d) ~\Big|~ i_s \in I_{t_s}, s = 1,...,d \Big\}, \\
  \bXi_\bt & = & \Big\{ (\xi_1,...,\xi_d) ~\Big|~ \xi_s \in \Xi_{t_s}, s = 1,...,d \Big\}
\end{eqnarray*}
and the tensor product functions
\begin{eqnarray*}
 v_{\bt,\bi}^\slin(\bx) & := & \prod_{s=1}^d v_{t_s,i_s}^\slin(x_s), \quad \bi \in \bI_\bt, \\
\varphi_{\bt,\bi}(\bx) & := & \prod_{s=1}^d \varphi_{t_s,i_s}(x_s), \quad \bi \in \bXi, \\
\end{eqnarray*}
where $\bx = ( x_1,...,x_d)$.
These constructions allow us to define the tensor product vector spaces (see Figure \ref{Picspaces})
\begin{eqnarray*}
 V_{n,d}^\sFull & := & \sspan \big\{ v_{\bt,\bi}^\slin ~\big|~ \max(\bt) \leq n, \bi \in \bXi_\bt  \big\}, \\
 V_\bt & := & \sspan \big\{ v_{\bt',\bi}^\slin ~\big|~ \bt' \leq \bt, \bi \in \bXi_{\bt'}  \big\}, \\
 W_\bt & := & \sspan \big\{  \varphi_{\bt,\bi} ~\big|~  \bi \in \bXi_{\bt}  \big\}, \\
 V_{\cD_n}^\slin & := & \sspan \big\{ v_{\bt,\bi}^\slin ~\big|~ |\bt| \leq n, \bi \in \bXi_\bt  \big\}, \\
 V_{\cD_n}^\sPre & := & \sspan \big\{ \varphi_{\bt,\bi} ~\big|~ |\bt| \leq n, \bi \in \bXi_\bt  \big\}. 
\end{eqnarray*}
Obviously, this results in $W_\bt \subset V_\bt$.

\begin{figure}
  \centering
  \includegraphics[width=0.25\textwidth]{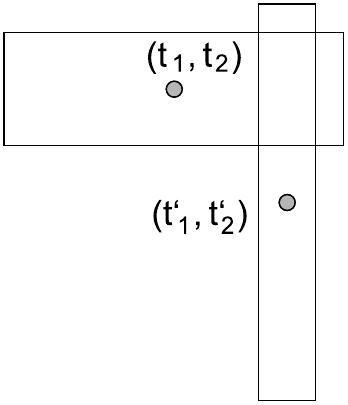}
  \caption{Example of the support of two basis functions that satisfy $|\max(\bt,\bt')| >n$ 
    and $|\max(t_1,t_2)| = | \bt | \leq n$ $|\max(t'_1,t'_2)| = | \bt' | \leq n$.}
\label{FuncSemiOrtho}
\end{figure}

$V_{n,d}^\sFull$ is the well-known standard space of multilinear finite element functions which
can be written as
\[
V_{n,d}^\sFull  = \sspan \big\{ v_{(n,...,n),\bi}^\slin ~\big|~ \bi \in \bI_{(n,...,n)}  \big\}.
\]
The sparse grid spaces $V_{\cD_n}^\slin$ and $V_{\cD_n}^\sPre$ are equal (see \cite{PflaumHartmann}).
However, the adaptive versions of these spaces are not equal. For reasons of simplicity, only the non-adaptive case is dealt with. In case of smooth functions, sparse and full grid have similar 
approximation properties
\begin{eqnarray*}
 \min_{v \in V_{n,d}^\sFull} \| u - v \|_{H^1} & \leq & C 2^{-n} \| u \|_{H^2} \quad \mbox{and} \\
 \min_{v \in V_{\cD_n}^\sPre} \| u - v \|_{H^1} & \leq & C \, n \, 2^{-n} 
 \left\| \frac{\partial^{2d} u}{\partial x_1^2 ... \partial x_d^2} \right\|_{L^2} ,
\end{eqnarray*}
where $C$ is a constant independent of $n$ and $u$.
However, the dimension of the corresponding spaces is completely different:
\begin{eqnarray*}
\mbox{dim} (V_{n,d}^\sFull) = O(2^{nd}) \quad \mbox{and} \quad 
\mbox{dim} (V_{\cD_n}^\sPre) = O( n^{d-1} 2^{n}) .
\end{eqnarray*}
Therefore, the aim of this paper is to find an efficient Galerkin discretization of Problem (\ref{PropOrg})
using the sparse grid space $V_{\cD_n}^\sPre$. To this end, the following lemma is an important observation:
\begin{lemma}[Semi-Orthogonality Property]
\label{SemiOrthoOne}
Let $\kappa$ be constant and $A = \mbox{diag}(\alpha_1,...,\alpha_d)$ a constant diagonal matrix.
Then for all indices $\bt$, $\bt'$, $\bi \in \bXi_\bt$, and $\bi' \in \bXi_{\bt'}$ such that
\begin{equation}
\label{assTheoOne}
 | \max(\bt,\bt') | > n ~ \mbox{and} ~ |\bt| \leq n, ~~|\bt'| \leq n 
\end{equation}
the following equation holds
 \[
  \int_\Omega (\nabla \varphi_{\bt,\bi})^T A \nabla \varphi_{\bt',\bi'} + 
           \kappa \varphi_{\bt,\bi} \varphi_{\bt',\bi'} ~d\bx = 0 .
\]
\end{lemma}
The consequence of this lemma is that pre-wavelet basis functions with overlapping support are orthogonal 
to each other (see Figure \ref{FuncSemiOrtho}). This orthogonality property is
the motivation of the following discretization:

\begin{discretization}[Semi-Orthogonality]
\label{DisB} 
\noindent
Let $f \in L^2(\Omega)$ and $\kappa \in L^\infty(\Omega), \kappa \geq 0$ be given. Then, let us define
\[
 a(u,v) := \int_\Omega (\nabla u)^T A(\bx) \nabla v + \kappa(\bx) u v ~d\bx 
\]
and
\begin{eqnarray*}
a_n^{\semiortho} : V_{\cD_n}^\sPre \times V_{\cD_n}^\sPre & \to & \real \\ 
 a_n^{\semiortho}(\varphi_{\bt,\bi}, \varphi_{\bt',\bi'}) & := & 
\left\{ 
\begin{array}{ccl}
a(\varphi_{\bt,\bi}, \varphi_{\bt',\bi'}) & \mbox{if} & | \max(\bt,\bt') | \leq n \\
0      & \mbox{if} & |\max(\bt,\bt')| > n
\end{array}
\right. .
\end{eqnarray*}
Find $u_{\cD_n}^\sPre \in V_{\cD_n}^\sPre$ such that
\begin{equation}
\label{equDisSemi}
 a_n^{\semiortho} ( u_{\cD_n}^\sPre , v_h ) = \int_\Omega f v_h ~d\bx \quad \forall v_h \in V_{\cD_n}^\sPre.
\end{equation}
\end{discretization}
In \cite{PflaumHartmann} we analyzed the convergence of this discretization for the Helmholtz problem with
variable coefficients with respect to the $H^1$-norm. This paper shows how to obtain an efficient
algorithm for solving the corresponding linear equation.
\section{Basic Notation}
\vspace{-2pt}
\label{Basic}
The difficulty in explaining sparse grid algorithms is that the matrices in these algorithms are applied to
vectors with varying sizes. Thus, describing these matrices in a mathematical correct form
leads to a non-trivial notation. A second problem appears in the case of adaptive grids. All sparse grid 
algorithms have a recursive structure using a tree data structure. Explaining such algorithms in a mathematical
and clear notation is difficult. Therefore, we restrict ourselves to non adaptive sparse 
grids and assume that a sparse grid is a union of  semi-coarsened full grids.

Furthermore, we introduce a notation that is based on operators on vector
spaces and its dual space. Assume that the finite element solution $u_\bt \in V_\bt$ ids searched such that
\[ 
 a_n^{\semiortho} ( u_\bt , v_\bt ) = \int_\Omega f v_\bt ~d\bx \quad \forall v_\bt \in V_\bt .
\]
Then $u_\bt$ is contained in the vector space $V_\bt$, but the mappings
\[
 v \mapsto \int_\Omega f v ~d\bx, \quad v \in V_\bt
\]
and
\[
 v \mapsto a_n^{\semiortho} ( w , v ), \quad v \in V_\bt
\]
are contained in the dual space $V'_\bt$.
To store elements of $V_\bt$ or
functionals of $V'_\bt$ as a vector, assume that
$n_\bt$ is the number of grid points on level $\bt$
\[
 n_\bt = \dim ( V_\bt )
\]
Moreover, assume that the data of a sparse grid algorithm is stored on various suitable full grids.
A corresponding global array is:
\begin{eqnarray*}
 \vec{U}_n & := & \left( U_\bt \right)_{|\bt| \leq n} \\
 U_\bt & \in  & \real^{n_\bt}
\end{eqnarray*}
The vector $U_\bt$ is used to store data of different mathematical objects.
One possibility is to describe a function in $V_\bt$ by the vector $U_\bt$. Another possibility
is to describe a functional of $V'_\bt$ by $U_\bt$. The notation for a corresponding assignment operator is:
\[
\left( U_\bt  \setto u \right) \;\; \text{for} \; u \in V_\bt
\]
and
\[
\left( U_\bt  \setto  f \right)  \;\; \text{for} \; f \in V_\bt'
\]
Some examples are: \\

\noindent
\quad {\bf Example 1} Let $u = \sum_{\bi \in \bI_\bt } c_{\bt,\bi} v_{\bt,\bi}^\slin \in V_\bt $. 
Then,
\[
 \left( U_\bt  \setto u \right)  
~~~ :\Leftrightarrow  ~~~  U_\bt = (c_{\bt,\bi})_{\bi \in \bI_\bt} .
\]
\quad {\bf Example 2} Let $u = \sum_{\bi \in \bXi_\bt } c_{\bt,\bi} \varphi_{\bt,\bi}  \in W_\bt \subset V_\bt $. 
Then, we write
\[
 \left( U_\bt  \setto u \right)  
~~~ :\Leftrightarrow  ~~~  U_\bt = \left(  \begin{array}{ll} c_{\bt,\bi} & \mbox{if $\bi \in \bXi_\bt$} \\
                                        0           & \mbox{else}.
                     \end{array} \right)   .
\]
\quad {\bf Example 3} Let $f \in V_\bt'$ Then, we write
  \[
 \left( U_\bt  \setto f \right)  
~~~ :\Leftrightarrow  ~~~  U_\bt = (f(v_{\bt,\bi}))_{\bi \in \bI_\bt} .
\] 

For describing our algorithms, we introduce a special operator $\cB$, which will be called
{\em back construction operator}. This operator reconstructs the mathematical object $u$
which was used to set values in a vector $U_\bt$. This implies the following property
of the back construction operator $\cB$:
\[
\left( U_\bt  \setto u \right)  ~~~ \Rightarrow ~~~
 u = \cB \left( U_\bt \right)  .
\]
Therefore, if $U_\bt$ was set by $u$ in an assignment $\left( U_\bt  \setto u \right)$ during
execution of an algorithm, then $\cB$  reconstructs $u$ in a later execution of this algorithm. 
Using this notation, one can describe algorithms without knowing 
how information on $u$ were stored in $U_\bt$.
\section{Basic Operators}
\vspace{-2pt}

\label{BasicOp}  
The algorithms in this paper are mainly based on well-known one-dimensional operators.
Briefly recall these operators:

\noindent
{\em 1. Prolongation} \\
A prolongation in direction $s$ can be described by
\[
  W_{\bt} \setto  I_\bt \left( \cB(W_{\bt - e_s}) \right) ,
\]
where $\cB(W_{\bt - e_s}) \in V_{\bt - e_s}$ is given and $e_s$ is the unit vector in direction $s$.
This operator is described in matrix form
\begin{eqnarray*}
 ( c^\coa_{\bi} )_{\bi \in \bI_{\bt-e_s}} \in \real^{|\bI_{\bt-e_s}|} , & ~~~~ & 
  \cB(W_{\bt-e_s}) = \sum_{\bi \in \bI_{\bt -e_s}}  c^\coa_{\bi} v_{\bt-e_s,\bi},  \\
  ( c_{\bi} )_{\bi \in \bI_\bt} \in \real^{|\bI_\bt|} , & ~~~~ & 
  \cB(W_\bt) = \sum_{\bi \in \bI_{\bt}}  c_{\bi} v_{\bt,\bi} .
\end{eqnarray*}
The prolongation in direction $s$ acts on these vectors as follows
\begin{eqnarray}
\label{prolMat}
  ( c_\bi )_{\bi \in \bI_{\bt}} & = &  M^\prol_s ( c^\coa_\bi )_{\bi \in \bI_{\bt-e_s}}, 
  \quad \quad M_s^\prol = \otimes_{i=1}^{s-1} \Id \otimes M^\prol \otimes \otimes_{i=s+1}^{d} \Id,
\end{eqnarray}
where $\Id$ is the identity matrix and $M^\prol$ is the matrix
\[
 M^\prol = \left(
\begin{array}{ccccccccccccc}
  1 &  &   \\[1ex]
  \frac{1}{2} & \frac{1}{2} & \\[1ex]
              &   1         & \\
              &      &    \ddots     \\
              &      & &  1 &  & \\[1ex]
              &      & &  \frac{1}{2} & \frac{1}{2}  \\
              &      & &              &       1               
\end{array} \right).
\]

\noindent
{\em 2. Restriction of the right hand side} \\
A restriction of the right hand side in direction $s$ can be described by
\[
  G_{\bt-e_s} \setto  \big( v \mapsto \cB(G_{\bt}) (v), ~ v \in V_{\bt-e_s} \big) ,
\]
where $\cB(G_\bt) \in V'_\bt$ is given.
To describe the matrix form of this operator let
\begin{eqnarray*}
 ( g_\bi )_{\bi \in \bI_\bt} \in \real^{|\bI_\bt|} , & ~~~~ & g_\bi = \cB(G_\bt)(v_{\bt,\bi}), \bi \in \bI_\bt, \\
 ( g^\coa_\bi )_{\bi \in \bI_{\bt-e_s}} \in \real^{|\bI_{\bt-e_s}|} , & ~~~~ & g_\bi = \cB(G_{\bt-e_s})(v_{{\bt-e_s},\bi}), \bi \in \bI_{\bt-e_s},  
\end{eqnarray*}
The restriction in direction $s$ acts on these vectors as follows
\begin{eqnarray}
\label{resMat}
  ( g^\coa_\bi )_{\bi \in \bI_{\bt-e_s}} & = &  M^\res_s ( g_\bi )_{\bi \in \bI_\bt}, 
  \quad \quad M_s^\res = \otimes_{i=1}^{s-1} \Id \otimes M^\res \otimes \otimes_{i=s+1}^{d} \Id,
\end{eqnarray}
where $\Id$ is the identity matrix and $M^\res$ is the matrix
\[
 M^\res = \left(
\begin{array}{ccccccccccccc}
  \frac{1}{2} & 1 &  \frac{1}{2} \\[1ex]
           &      &        &  \ddots     \\
               & &&             &  \frac{1}{2} & 1 & \frac{1}{2}  \\[1ex]
   &&                &&&     &  \frac{1}{2} & 1 & \frac{1}{2} & \\
   &&         &&&&     &    &      & \ddots \\              
        & &    & & &&& &  &  & \frac{1}{2} & 1 & \frac{1}{2}             
\end{array} \right).
\]

\noindent
{\em 3. Transformation to pre-wavelets} \\
Let $s$ be a direction  such that $1 \leq s \leq d$. Furthermore, assume that $U_\bt$ is given
such that $\cB(U_\bt) \in V_\bt$ and
\begin{itemize}
 \item $U_\bt$ stores the coefficients of $\cB(U_\bt)$ in pre-wavelet form or nodal basis form in the directions
$\tilde{d} \not= s$, and
\item $U_\bt$ stores the coefficients of $\cB(U_\bt)$ in nodal basis form in the direction $s$.
\end{itemize}
Now assume that the pre-wavelet coefficients are calculated in direction $s$ with respect to
level $t_s$ of  $\cB(U_\bt)$  and the resulting vector is stored in $H_\bt$. 
The corresponding assignment can be written as
\begin{equation}
\label{AssA}
   H_\bt \setto Q^s_{\bt} (\cB(U_\bt) )  .
\end{equation}
Here $Q^s_{\bt}$ is the $L^2$ projection operator onto the space
\[
 \otimes_{i=1}^{s-1} V_{t_i} \otimes W_{t_s} \otimes \otimes_{i=s+1}^{d} V_{t_i}.
\]
In matrix form, assignment (\ref{AssA}) can be written as
\[
   \otimes_{i=1}^{s-1} \Id \otimes M^\sPre \otimes \otimes_{i=s+1}^{d} \Id,
\]
where $\Id$ are suitable identity matrices and $M^\sPre$ is the matrix of the 1-dimensional case.
To describe the matrix $M^\sPre$ of the 1-dimensional case, let
\begin{eqnarray*}
 ( u_i )_{i \in I_t} \in \real^{|I_t|} , & ~~~~ & u_i = \cB(U_t) = \sum_{i \in I_t} v_{t,i},  \\
 ( h_i )_{i \in \Xi_t} \in \real^{|\Xi_t|} , & ~~~~ & h_i = \cB(H_t) = \sum_{i \in \Xi_t} \varphi_{t,i} ,
\end{eqnarray*}
where $t = t_s$.
The matrix $M^\sPre$ performs the following mapping
\[
 ( h_i )_{i \in \Xi_t} =  M^\sPre ( u_i )_{i \in I_t} .
\]
Let $\cR^\sPre$ be the following restriction operator which takes only the pre-wavelet coefficients:
\begin{eqnarray*}
 \cR^\sPre \left( ( h_i )_{i \in I_t} \right) & := &  ( h_i )_{i \in \Xi_t}.
\end{eqnarray*}
Then
\[
 M^\sPre = \cR^\sPre M^{-1},
\]
where $M$ is the matrix
\begin{equation}
\label{matrixM}
 M = \left(
\begin{array}{ccccccccccccccc}
  \frac{9}{10} & \frac{1}{2} &  \frac{1}{10} \\[1ex]
 -\frac{3}{5}  & 1           & -\frac{3}{5}  &  \\[1ex]
  \frac{1}{10} & \frac{1}{2} &   1            &  \ddots & \\
               &             &  -\frac{3}{5}  & \ddots  & \frac{1}{2}  &  \frac{1}{10}  \\
               &             &  \frac{1}{10}  &  \ddots & 1            &  -\frac{3}{5}  \\
               &             &                &         &  \frac{1}{2} & \frac{9}{10}   \\[1ex]
\end{array} \right).
\end{equation}
This means that the assignment (\ref{AssA}) has to be implemented by inverting the matrix
(\ref{matrixM}) in direction $s$ and taking only the resulting pre-wavelet coefficients.

\noindent
{\em 4. Discretization stencil} \\
Let $U_\bt$ be a given vector on a semi-coarsened full subgrid such that $\cB(U_\bt) \in V_\bt$.
Then an assignment involving the bilinear form of the operator is:
\[
Z_{\bt}  \setto  \big( v \mapsto a( \cB(U_\bt),v), ~ v \in V_{\bt} \big) .
\]
Obviously, the computation corresponding to this assignment requires a 9-point-stencil in the 2-dimensional case and
a 27-point-stencil in the 3-dimensional case. As an example, for meshsize $h=h_x=h_y$ in x-direction
and y-direction this stencil is
\[
\frac{1}{h^2} \left[\begin{array}{ccc}
  -1 & -1 & -1 \\
  -1 &  8 & -1 \\
  -1 & -1 & -1 \\
\end{array} \right] ,
\]
where 
\[
 a(u,v) = \int_\Omega \nabla u \nabla v ~d(x,y).
\]

\section{Calculation of pre-wavelet coefficients}
\vspace{-2pt}

Let $ \left( F_{\bt} \right)_{ | \bt \ \leq n  } $ be a function evaluated on the sparse grid of depth $ n $
such that $ F_{\bt} = \left( f \left( x_{ \bt,\bi } \right) \right)_{\bi \in \bI_{\bt}} $ where $ x_{\bt,\bi} $
is the grid point on level $ \bt $ with index $ \bi $.
Algorithm \ref{prew_coeff} calculates the pre-wavelet decomposition with coefficients
$ \left( C_{\bt} \right)_{ | \bt | \leq n} =
\left( c_{\bt,\bi} \right)_{ | \bt | \leq n, \bi \in \bI_{\bt}} $ given by
\begin{equation}
 f \left( x_{ \bt,\bi } \right) = \sum\limits_{| \bt | \leq n, \bi \in \bI_{\bt}} { c_{\bt,\bi} \varphi_{\bt,\bi}\left( x_{ \bt,\bi } \right) }
\end{equation}
for every sparse grid point $x_{ \bt,\bi } \in \cD_n$, where $\cD_n$ is the sparse grid of depth $n$.

The calculation of the wavelet coefficients works recursively through all the dimensions.
For this purpose, the coefficients of the highest dimension are calculated first, starting from the 
grid with maximum depth down to the coarsest grid.
After this, the algorithm continues with lower dimensions.

The calculation of the pre-wavelet coefficients on depth $ t $ in one dimension requires to solve a system of 
linear equations with $ 2 ^ {t + 1} +1 $ unknowns.
The inverse matrix $ M^{-1} $ can efficiently be computed using an LU decomposition since $M$ is a band matrix with bandwidth 3  (see Equation (\ref{matrixM})).


Now, let us explain the basic idea 
of the algorithm for dimension $ 0 \leq \tilde{d} < d  $ on level $ t_{\tilde{d}} $.
First, it calculates the pre-wavelet decomposition.
Then, the local hierarchical surplus is subtracted from all low order levels $ \tilde{t}_{\tilde{d}} < t_{\tilde{d}} $.
In the case of sparse grids, this requires the interpolation of the the hierarchical surplus for grids
which have no direct predecessor. To this end, 
the combination technique is applied to all direct neighbors in directions with the dimension 
larger than $ \tilde{d} $.

The same algorithm is used in reverse order to calculate a point-wise evaluation 
$ F_{\bt} = \left( f \left( x_{ \bt,\bi } \right) \right)_{\bi \in \bI_{\bt}} $ for a given set of 
pre-wavelet coefficients $ \left( C_{\bt} \right)_{ | \bt | \leq n} $.
Algorithm \ref{prew_back} starts on the coarsest level of the smallest dimension and accumulates the 
surpluses over all dimensions.

Observe that in Algorithm \ref{prew_back} and Algorithm \ref{prew_coeff} $\cU_{\bar{d}}(\bt)$ is
the upper part of $\bt$:
\[
 \cU_{\bar{d}}(\bt)  = (t_{\bar{d}+1}, t_{\bar{d}+2}, ... ,  t_d) .
\]

\begin{algorithm}[Calculate pre-wavelet decomposition]~
\label{prew_coeff}
\begin{tabbing}
 aa \= aa \=  aa \=  aa \=  aa \=  aa \kill
\Markiere{Input:} Let $\vec{U}_n = \left( U_\bt \right)_{|\bt| \leq n}$ be given in nodal format. \\
\> \Markiere{Call}  {\sc Pre-wavelet Algorithm}  ( $\vec{U}_n$, $d$). \\
\Markiere{Output:}  $\vec{U}_n$ in pre-wavelet format until dimension $\bar{d}$. \\
 \\
\Markiere{Function} {\sc Pre-wavelet Algorithm}
( $\vec{U}_m  :=  \left( U_\bt \right)_{|\bt|_{\bar{d}} \leq m, \cU_{\bar{d}}(\bt) =\bt_u}$, $\bar{d}$)  $\{$ \\
\> \Markiere{iterate} for $\bar{t}=m,...,0$ $\{$ \\
\> \> 1. Calculate pre-wavelet  coefficients in direction $\bar{d}$: \\
\> \>   \Markiere{for} every  $\bt$ with $t_{\bar{d}} = \bar{t}$ and $\cU_{\bar{d}}(\bt) =\bt_u$ do: \\
\> \>  \>   $ U_\bt \setto Q^{\bar{d}}_{\bt} (\cB(U_\bt) ) $ \\
\> \>  \>  $ W_\bt   =    U_\bt $ \\
\> \> 2. Subtract interpolated pre-wavelets on coarse grid \\
\> \>  \Markiere{iterate} for $t'=\bar{t},...,1$ $\{$ \\
\> \> \>  2.1 \Markiere{for} every  $\bt$ with $t_{\bar{d}} = t'$ and $\cU_{\bar{d}}(\bt) =\bt_u$ do: \\
\> \> \> \> $W_{\bt - e_{\bar{d}}} \setto I^{\bar{d}}_{\bt - e_{\bar{d}}}( \cB(W_\bt) )$   \\
\> \> \> \> $U_{\bt - e_{\bar{d}}} \setto \cB(U_{\bt - e_{\bar{d}}}) - ( \cB(W_{\bt - e_{\bar{d}}}) ) $ \\
\> \> \> 2.2 \Markiere{if} $\bar{d}>1$ then $\{$ \\
\> \> \> \>  \Markiere{for} every  $\bt$ with $|\bt|_{\bar{d}} = m$ and $t_{\bar{d}} = t'-1$ and  $\cU_{\bar{d}}(\bt) =\bt_u$ do: \\
\> \> \> \> \> 
     $ W_{\bt} \setto 
\sum_{ (\alpha_1,...,\alpha_{\bar{d}-1}) \in \{ 0,1 \}^{\bar{d}-1}  \wedge
      \atop (\alpha_{\bar{d}},...,\alpha_{d}) = (0,...,0) \wedge \alpha \not= (0,...,0) } 
  (-1)^{| \alpha | + 1}  I_{\bt} ( \cB(W_{\bt-\alpha}) ) $  \\
\> \> \> \> \> $ U_{\bt} \setto \cB(U_{\bt} ) - \cB( W_{\bt} ) $ \\
\> \> \> $\}$ \\
\> \>    $\}$ \\
\> \>  3. Recursion \\
\> \> \Markiere{if} $\bar{d}>1$ then $\{$ \\
\> \>  \> Define $m^{low} := m - \bar{t}$.  Define $\vec{U}^{low}_{m^{low}}  :=  \left( U_{\bt} \right)_{|\bt|_{\bar{d}-1} \leq m^{low}, t_{\bar{d}} = \bar{t},\cU_{\bar{d}}(\bt) =\bt_u}$ \\ 
\> \>   \> If $m^{low}>0$ \Markiere{call}  {\sc Pre-wavelet Algorithm}  ( $\vec{U}^{low}_{m^{low}}$, $\bar{d}-1$) \\
\> \>   $\}$ \\
 \>    $\}$ \\
$\}$ \\
End of Algorithm.  
\end{tabbing}
\end{algorithm}

\begin{algorithm}[Calculate back pre-wavelet decomposition]~
\label{prew_back}
\begin{tabbing}
 aa \= aa \=  aa \=  aa \=  aa \=  aa \kill
\Markiere{Input:} Let $\vec{U}_n = \left( U_\bt \right)_{|\bt| \leq n}$ be given in pre-wavelet format. \\
 \> \Markiere{Call}  {\sc Back Pre-wavelet Algorithm}  ( $\vec{U}_n$, $d$). \\
\Markiere{Output:}  $\vec{U}_n$ in nodal format until dimension $\bar{d}$. \\
\\
\Markiere{Function} {\sc Back Pre-wavelet Algorithm}
( $\vec{U}_m  :=  \left( U_\bt \right)_{|\bt|_{\bar{d}} \leq m, \cU_{\bar{d}}(\bt) =\bt_u}$, $\bar{d}$) $\{$ \\
\> \Markiere{iterate} for $\bar{t}=0,...,m$  $\{$ \\
\>  \>  1. Recursion \\
\>  \>  \Markiere{if} $\bar{d}>1$ then $\{$ \\
\> \>  \> Define $m^{low} := m - \bar{t}$.  Define $\vec{U}^{low}_{m^{low}}  :=  \left( U_{\bt} \right)_{|\bt|_{\bar{d}-1} \leq m^{low}, t_{\bar{d}} = \bar{t},\cU_{\bar{d}}(\bt) =\bt_u}$ \\ 
\>  \> \>  If $m^{low}>0$ \Markiere{call}  {\sc Back Pre-wavelet Algorithm}  ( $\vec{U}^{low}_{m^{low}}$, $\bar{d}-1$) \\
\>  \>   $\}$ \\
\>  \>  2. Transform back in direction $\bar{d}$: \\
\>  \>  \Markiere{if} $\bar{t} > 0$   then $\{$ \\ 
\>  \>  \> \Markiere{for} every  $\bt$ with $t_{\bar{d}} = \bar{t}$ and $\cU_{\bar{d}}(\bt) =\bt_u$ do: \\
\>  \> \> \>   $ W_\bt  \setto  I_{\bt}(\cB(U_\bt)) $ \\
\>  \> \> \>  $ U_\bt  \setto  \cB(W_\bt)  + I_{\bt}^{\bar{d}}(\cB(U_{\bt - e_{\bar{d}}})) $ \\
\>  \>  $\}$ \\
\> \> 3. Add interpolated pre-wavelets on coarse grid \\
\>  \>  \Markiere{iterate} for $t'=\bar{t},...,1$ $\{$ \\
\>   \> \>   3.1. \Markiere{for} every  $\bt$ with $t_{\bar{d}} = t'$ and $\cU_{\bar{d}}(\bt) =\bt_u$ do: \\
\>   \> \> \>   $W_{\bt - e_{\bar{d}}} \setto I^{\bar{d}}_{\bt - e_{\bar{d}}}( \cB(W_\bt) )$   \\
\>   \> \> \>   $U_{\bt - e_{\bar{d}}} \setto \cB(U_{\bt - e_{\bar{d}}}) + ( \cB(W_{\bt - e_{\bar{d}}}) ) $ \\
\>   \> \>   3.2. \Markiere{if} $\bar{d}>1$ then $\{$  \\
\>   \> \> \>  \Markiere{for} every  $\bt$ with $|\bt|_{\bar{d}} = m$ and $t_{\bar{d}} = t'-1$ and  $\cU_{\bar{d}}(\bt) =\bt_u$ do: \\
\>   \> \> \> \> $  W_{\bt} \setto  \sum_{ (\alpha_1,...,\alpha_{\bar{d}-1}) \in \{ 0,1 \}^{\bar{d}-1}  \wedge
                                  \atop (\alpha_{\bar{d}},...,\alpha_{d}) = (0,...,0) \wedge \alpha \not= (0,...,0)     } 
(-1)^{| \alpha | + 1}  I_{\bt} ( \cB(W_{\bt-\alpha}) ) $ \\  
\>  \> \> \> \> $U_{\bt}  \setto  \cB(U_{\bt} ) + \cB( W_{\bt} ) $ \\
\> \>  \>  $\}$ \\ 
\> \>      $\}$ \\
\>  $\}$ \\
$\}$ \\
End of Algorithm.  
\end{tabbing}
\end{algorithm}
\section{Matrix multiplication}
\vspace{-2pt}

\label{SectionMatrixVector}
Let
\begin{equation}
 \label{StiffnessMat}
\cA_n = \Big( a_n^{\semiortho} (\varphi_{\bt,\bi},\varphi_{\bt',\bi'}) 
\Big)_{ {|\bt| \leq n, \bi \in \bXi_\bt} \atop {|\bt'| \leq n, \bi' \in \bXi_{\bt'}} }.
\end{equation}
be the stiffness matrix of Discretization \ref{DisB}.
The main difficulty is to construct an algorithm that efficiently evaluates a matrix vector multiplication
with matrix $\cA_n$. To construct such an algorithm, the notation introduced in Section \ref{Basic} is applied.
The following recursive algorithm is obtained:

\begin{algorithm}[Sparse Grid Matrix Multiplication]~
\label{AlgMatrix}
\begin{tabbing}
 aa \= aa \=  aa \=  aa \=  aaaaaa \= aaaaaaaa \=  aa \kill
\Markiere{Input:} Let $\vec{U}_n = \left( U_\bt \right)_{|\bt| \leq n}$ be given in pre-wavelet format. \\
 \>    $F_{\bt} \setto  0  ~~ \mbox{$\forall \bt$ }$ \\
 \> \Markiere{for} every $R \subset \{ 1,2,...,d \}$ do $\{$ \\
 \> \> Reorder the directions $1,...,d$ such that  $R = \{ 1,2,...,D \}, \quad 0 \leq D \leq d$. \\
 \> \>   $\vec{H}_n : = \left( H_\bt \right)_{|\bt| \leq n} = \vec{U}_n$ \\   
 \> \> \Markiere{Call}  {\sc Sparse Grid Matrix Multiplication}  ( $\vec{H}_n$, $d$,$D$). \\
 \> $\}$ \\
\Markiere{Output:}  $\vec{F}_n$ in complete pre-wavelet format. \\ 
\Markiere{Function} {\sc Sparse Grid Matrix Multiplication} 
( $\vec{H}_m  :=  \left( H_\bt \right)_{|\bt|_{\bar{d}} \leq m, \cU_{\bar{d}}(\bt) =\bt_u}$, $\bar{d}$ , $D$) \\
\Markiere{if} $(\bar{d}==0)$ $\{$ // 1. Deepest point of recursion \\
\>   \Markiere{if} $(D==d)$  do \> \> \> \> $ Z_{\bt}  \setto  \big( v \mapsto a( \cB(H_\bt),v), ~ v \in V_{\bt} \big) $ \\
\>   \Markiere{else} \> \> \> \> $ Z_{\bt}  \setto  \cB(G_{\bt})$  \\
\>   $ F_{\bt}  \setto \cB(F_{\bt}) + \big( \varphi \mapsto \cB(Z_{\bt}) (\varphi), ~ \varphi \in W_\bt \big) $ \\
 $\}$ \Markiere{else if} $\bar{d} > D$ $\{$ // 2. Interpolate from coarser grids   \\
\>  2.1. Recursive call on coarsest grid. \\
\>   \> Define $\vec{H}^{low}_{m}  :=  \left( H_\bt \right)_{|\bt|_{\bar{d}-1} \leq m, t_{\bar{d}} = 0,\cU_{\bar{d}}(\bt) =\bt_u}$ \\
\>   \> \Markiere{call} 
   {\sc Sparse Grid Matrix Multiplication}  ( $\vec{H}^{low}_{m}$, $\bar{d}-1$, $D$) \\
\>  \Markiere{iterate} for $\bar{t}=1,...,m$  $\{$ \\
\>  \>  2.2. Prolongate in direction $\bar{d}$ \\
\>    \> \Markiere{for} every  $\bt$ with $t_{\bar{d}} = \bar{t}$ and $\cU_{\bar{d}}(\bt) =\bt_u$ do  \\
 \>   \> \>  $ H_\bt  \setto I_{\bt}^{\bar{d}}(\cB(H_\bt))  + I_{\bt}^{\bar{d}}(\cB(H_{\bt - e_{\bar{d}}})) $ \\
\>  \>   2.3. Recursive call on finer grids. \\
 \> \> Define $m^{low} := m - \bar{t}$.  Define $\vec{H}^{low}_{m^{low}}  :=  \left( H_\bt \right)_{|\bt|_{\bar{d}-1} \leq m, t_{\bar{d}} = \bar{t},\cU_{\bar{d}}(\bt) =\bt_u}$ \\
\>   \>  \Markiere{if} $m^{low}>0$ \Markiere{call} 
   {\sc Sparse Grid Matrix Multiplication}  ( $\vec{H}^{low}_{m^{low}}$, $\bar{d}-1$, $D$) \\
 \>  $\}$ \\
$\}$ \Markiere{else if} $\bar{d} \leq D$ $\{$ // 3. Restrict fine pre-wavelet functions  \\
\>  3.1. Recursive call on finest grid. \\
\>  \>    $G_{\bt} \setto  0$  ~~ $\forall \bt$, $t_{\bar{d}} = m,\cU_{\bar{d}}(\bt) =\bt_u$ \\
\>   \> Define $\vec{H}^{low}_{0}  :=  \left( H_\bt \right)_{|\bt|_{\bar{d}-1} = 0, t_{\bar{d}} = m,\cU_{\bar{d}}(\bt) =\bt_u}$ \\
\>   \> \Markiere{call} 
   {\sc Sparse Grid Matrix Multiplication}  ( $\vec{H}^{low}_{0}$, $\bar{d}-1$, $D$) \\
\>  \Markiere{iterate} for $\bar{t}=m-1,...,0$  $\{$ \\
\>  \>   3.2. Restrict right hand side. \\
\>   \>  \Markiere{for} every  $\bt$ with $t_{\bar{d}} = \bar{t}+1$ and $\cU_{\bar{d}}(\bt) =\bt_u$ do: \\
\>   \>  \>  $Q_{\bt} \setto  I_\bt ( \cB(H_\bt) ) $ \\
 \>    \> \>    $ G_{\bt - e_{\bar{d}}} \setto   \big( v \mapsto \cB(G_{\bt})(v) + a( \cB(Q_{\bt}),v), ~ v \in V_{\bt- e_{\bar{d}}} \big)$ \\
\>   \> Define $m^{low} := m - \bar{t}$.  Define $\vec{H}^{low}_{m^{low}}  :=  \left( H_\bt \right)_{|\bt|_{\bar{d}-1} \leq m^{low}, t_{\bar{d}} = \bar{t},\cU_{\bar{d}}(\bt) =\bt_u}$ \\
\>    \> \Markiere{if} $m^{low}>0$ \Markiere{call} 
   {\sc Sparse Grid Matrix Multiplication}  ( $\vec{H}^{low}_{m^{low}}$, $\bar{d}-1$, $D$) \\
$\}$ \>  $\}$ \\
End of Algorithm.  
\end{tabbing}
\end{algorithm}

\begin{figure}
  \subfigure[$ \left( pro,pro \right) $]{
    \label{YYY_a}\includegraphics[width=0.17\textwidth]{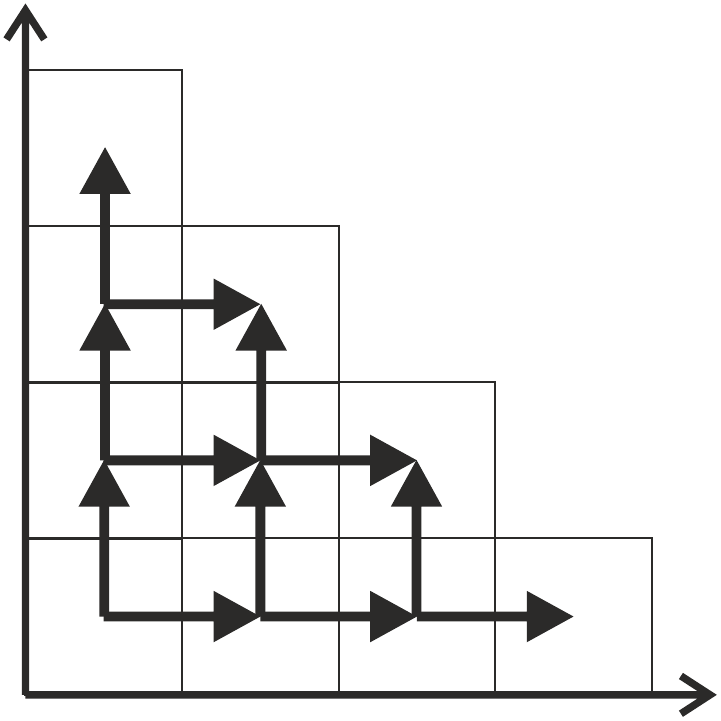}
  }
  \subfigure[$ \left( pro,res \right) $]{
    \label{YYY_b}\includegraphics[width=0.17\textwidth]{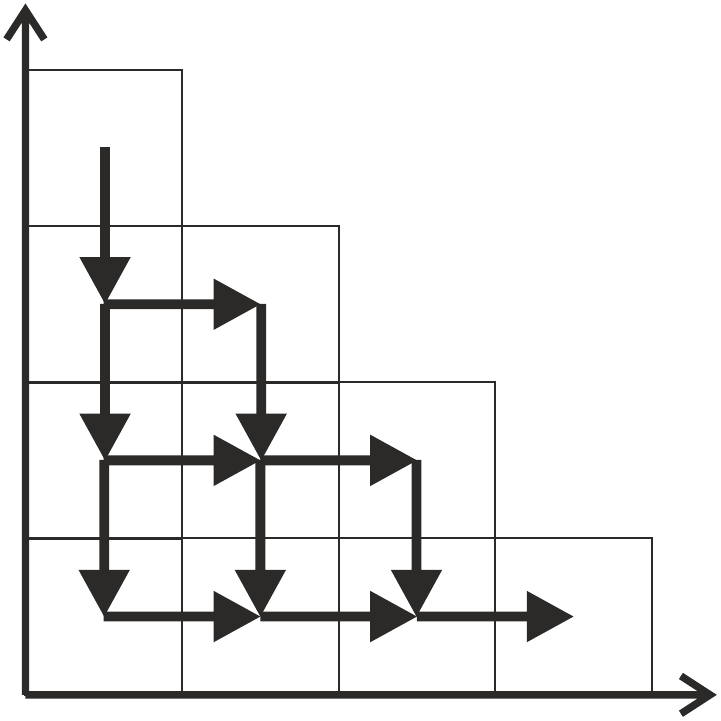}
  } 
  \subfigure[$ \left( res,pro \right) $]{
    \label{YYY_c}\includegraphics[width=0.17\textwidth]{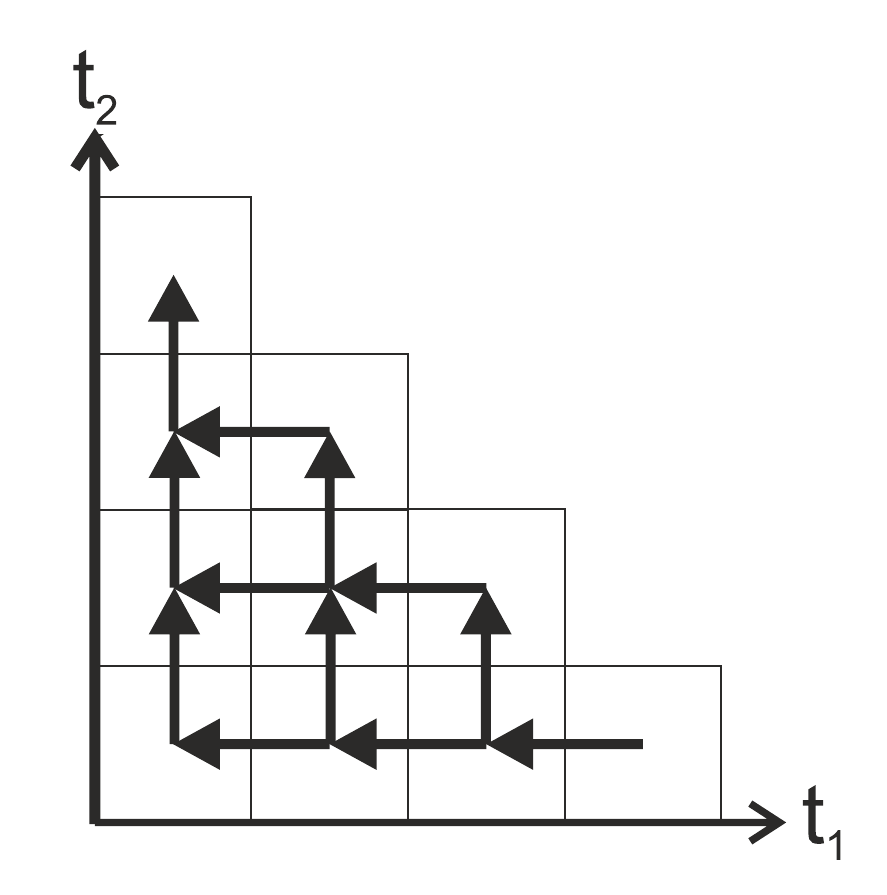}
  }
  \subfigure[$ \left( res,res \right) $]{
    \label{YYY_d}\includegraphics[width=0.17\textwidth]{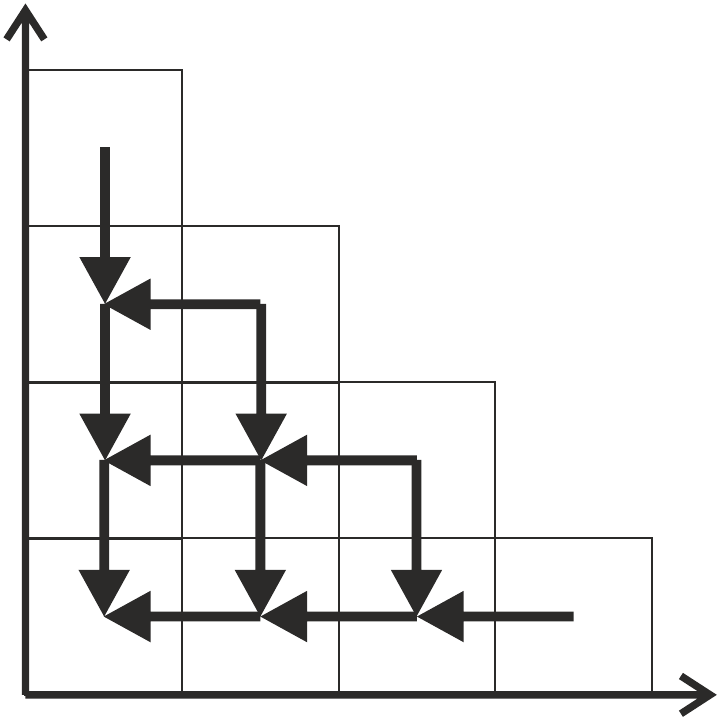}
  }
  \subfigure[Splitting]{
    \label{YYY_s}\includegraphics[width=0.17\textwidth]{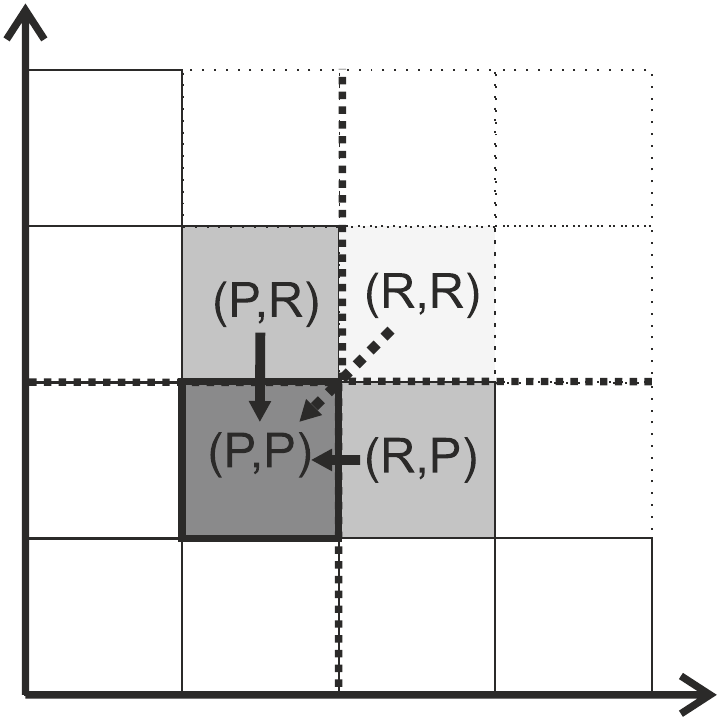}
  }  
  \caption{Matrix vector multiplication on sparse grids in the 2-dimensional case. 
   a) - d) Calculation of part $P=\{1,2\}$, $P=\{ 1 \}$, $P=\{ 2\}$, $P=\emptyset$.
   e) Summation of the result of the $4$ parts.}
  \label{TwoDMVMult}
\end{figure}

To explain the concept of this algorithm, start with a 1-dimensional observation.
Assume that the following pre-wavelet decomposition is given:
\[
 u = \sum_{t'=1}^n  w_{t'} \quad 
 \mbox{where} \quad w_{t'} = \sum_{i \in \Xi_{t'} }c_{t',i} \varphi_{t',i} .
\]
Then, we can split the matrix vector multiplication in two parts:
\[
 a(u,v_t) = \sum_{t' > t} a (w_{t'}, v_t)  ~~~  + ~~~ a ( \sum_{t' \leq t} w_{t'}, v_t) .
\]
Algorithm \ref{AlgMatrix} calculates these two parts separately. We call these two parts
{\em restriction part} and {\em prolongation part}.

The restriction part $\sum_{t' > t} a (w_{t'}, v_t) $ 
is calculated in the  function {\sc Sparse Grid Matrix Multiplication} with input parameter $D=1$.
The essential calculations are done in Step 3.2 and Step 1.
The recursive call of the algorithm applies the discretization stencil on $H_\bt$ (see Step 1) first
and then performs several restrictions (see Step 3.2).

The prolongation part $a ( \sum_{t' \leq t} w_{t'}, v_t)$
is calculated in the  function {\sc Sparse Grid Matrix Multiplication} with input parameter $D=0$.
Here the essential calculations are done in Step 2.2 and Step 1.
First $H_\bt$ is prolongated on finer grids (see Step 2.2) and then the discretization stencil is applied
(see Step 1).

The $d$-dimensional case is a combination of these two cases in all dimensional directions.
In each direction $1\leq \bar{d} \leq d$, we either can perform a restriction or prolongation.
All directions corresponding to restrictions are denoted by 
\[
R \subset \{ 1,...,d \}.
\]
Now, the recursive structure of Algorithm \ref{AlgMatrix} can be explained. Assume that $u$ is given in the following
form
\begin{eqnarray*}
 u & = & \sum_{|\bt'| \leq n} w_{\bt'}, \quad \quad \mbox{where} ~   w_{\bt'} 
         =   \sum_{\bi \in \Xi_\bt} \varphi_{\bt,\bi} c_{\bt,\bi}.
\end{eqnarray*}
This form corresponds to a decomposition of $u$ in its pre-wavelets.
Then a short calculation shows:
\begin{eqnarray} 
\label{sumCasesMat}
 a(u,\varphi_{\bt,\bi} ) & = & \sum_{R \subset  \{ 1,2,\ldots,d \} }  
      ~~~~~  \sum_{ s \in R }  \sum_{ t_s' > t_s}  
           a \left(  \sum_{ s \not\in R } \sum_{ t_s' \leq t_s} w_{\bt'}, \varphi_{\bt,\bi} \right).
\end{eqnarray}
In the 2-dimensional case this sum leads to $4$  parts:
\begin{eqnarray*}
 \lefteqn{ a(u,\varphi_{(t_1,t_2),(i_1,i_2)} ) = \sum_{ {t_1' > t_1} ,~  {t_2' > t_2}}  
           a \left(  w_{(t_1',t_2')}, \varphi_{(t_1,t_2),(i_1,i_2)} \right) } \\
                                     & + &  \sum_{ t_1' > t_1}  
           a \left(  \sum_{{t_2' \leq t_2}} w_{(t_1',t_2')}, \varphi_{(t_1,t_2),(i_1,i_2)} \right) 
                                      +   \sum_{ t_2' > t_2 }  
           a \left(  \sum_{{t_1' \leq t_1}} w_{(t_1',t_2')}, \varphi_{(t_1,t_2),(i_1,i_2)} \right) \\
  & + &    a \left(  \sum_{{t_1' \leq t_1} ,~ {t_2' \leq t_2}} w_{(t_1',t_2')}, \varphi_{(t_1,t_2),(i_1,i_2)} \right) .
\end{eqnarray*}
Figure \ref{TwoDMVMult} depicts the calculation of these four parts.

The d-dimensional case (\ref{sumCasesMat}) consists of $2^d$ cases. Each is calculated
by calling the function {\sc Sparse Grid Matrix Multiplication} with input parameter $d$ and $D$. Here $D$ is the number
of restriction directions. In order to simplify the notation, the restriction directions are reordered such that
\[
 R = \{ 1, \ldots,D \}.
\]
However, this is only a notational trick. An actual implementation of the algorithm should not perform
such a reordering.

\begin{figure}
 \centering
 \subfigure[2D]{
   \label{XXX_a}\includegraphics[width=2cm]{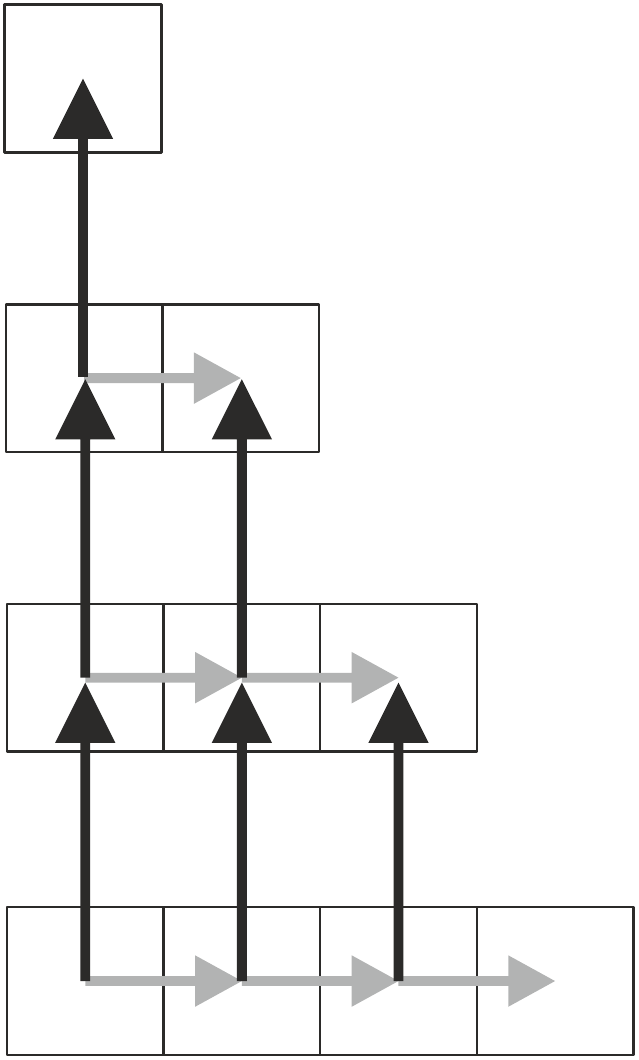}
 }
 \subfigure[3D]{
   \label{XXX_b}\includegraphics[width=2cm]{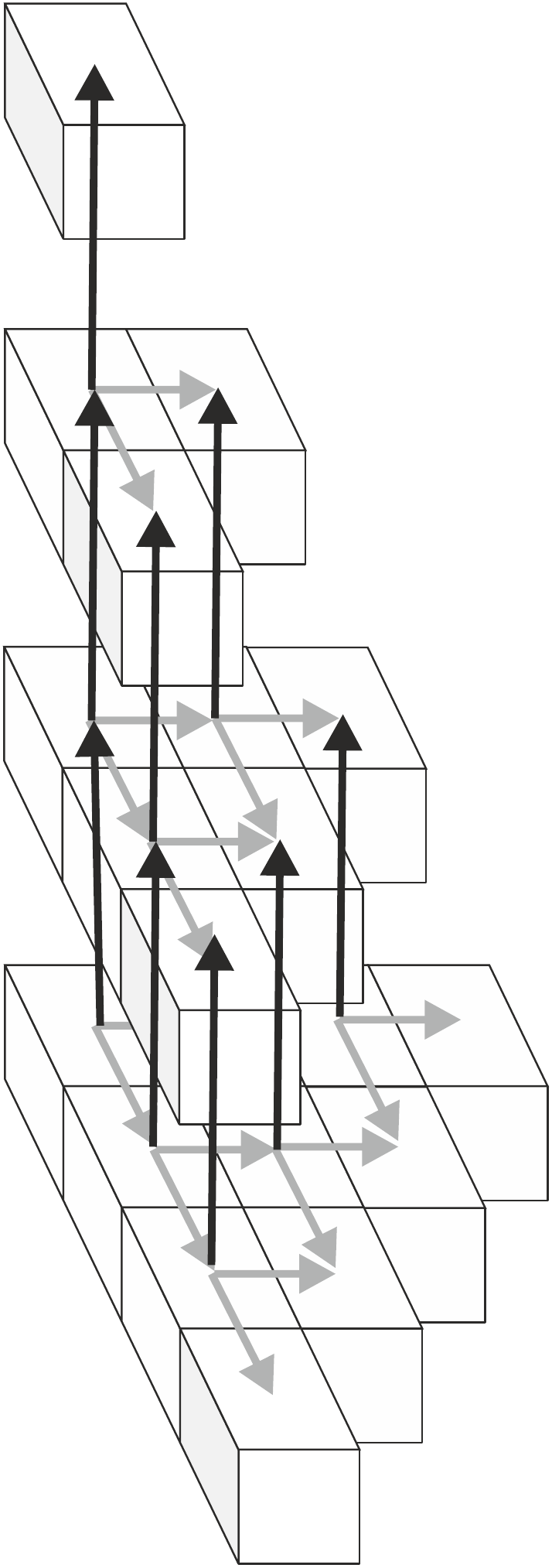}
 } 
 \caption{Prolongation step. \ref{XXX_a}: 2-dimensional case. 
  \ref{XXX_b}: 3-dimensional case.}
    \label{FigProll} 
\end{figure}

Observe that the recursive structure of Algorithm \ref{AlgMatrix} leads to the following sequential
calculations
\begin{enumerate}
 \item Prolongation in directions $D+1, \ldots, d$. This is depicted in Figure \ref{FigProll}.
 \item Application of discretization stencil on level $\bt$ in Step 1.
 \item Restriction in direction $1,\ldots , D$.
\end{enumerate}

\begin{figure}
  \center
  \subfigure[full grid]{
    \label{ZZZ_full}\includegraphics[width=0.3\textwidth]{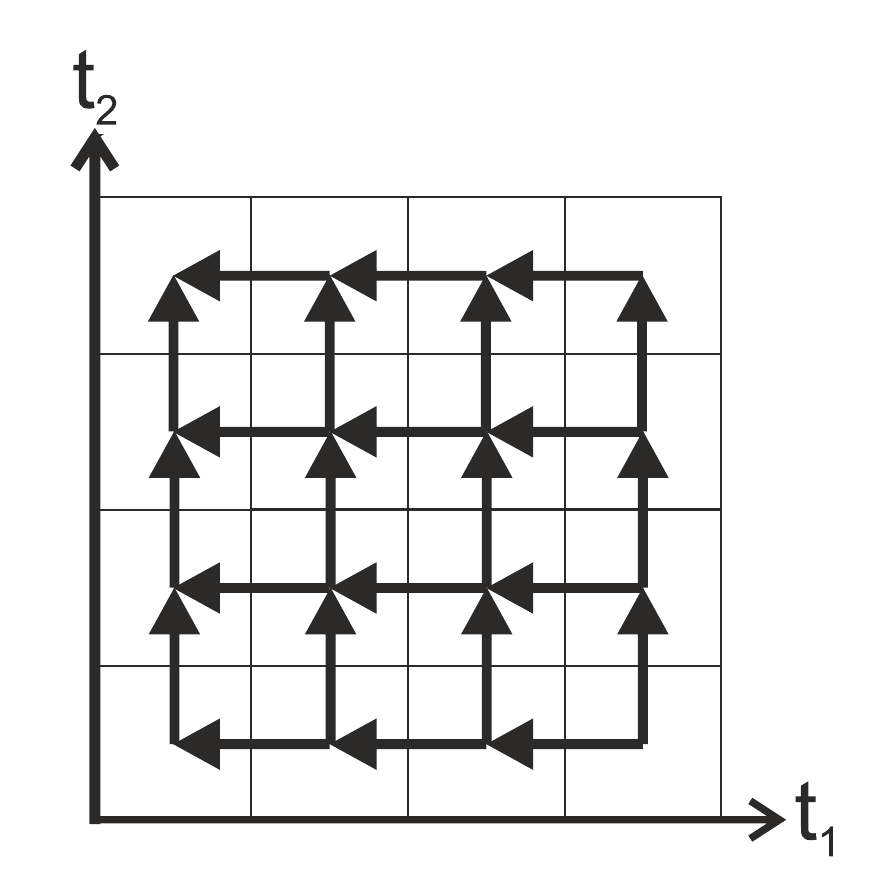}
  } 
  \subfigure[sparse grid]{
    \label{ZZZ_sparse}\includegraphics[width=0.3\textwidth]{./images/four_cases_RP_sparse.pdf}
  } \\
  \subfigure[missing path]{
    \label{ZZZ_missing}\includegraphics[width=0.3\textwidth]{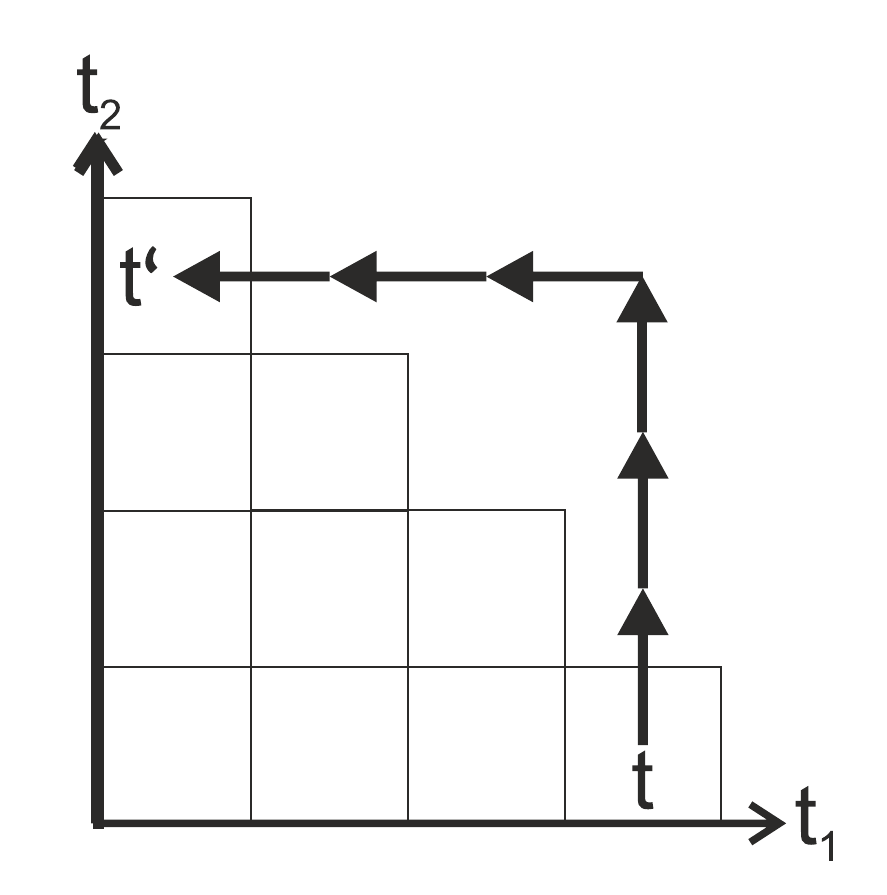}
  } 
  \subfigure[semi-orthogonality]{
    \label{ZZZ_semi}\includegraphics[width=0.3\textwidth]{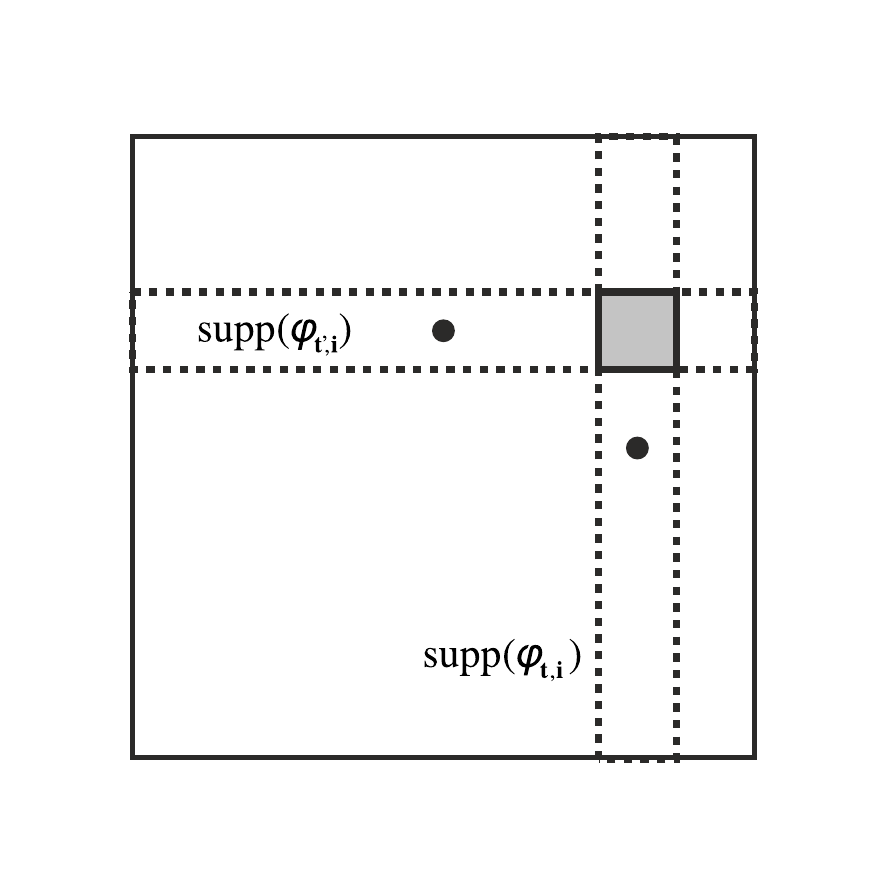}
  }  
  \caption{Prolongation-restriction in 2d on a sparse grid ignores the overlapping of $ \mathbf{t} $ and $ \mathbf{t}' $. Calculation would require a path over the full grid but can be ignored due to semi-orthogonality.}
  \label{four_cases}
\end{figure}

Now, two important question are: 
\begin{center}
{\em
 Does Algorithm \ref{AlgMatrix} correctly calculate a matrix vector multiplication? \\
 How is semi-orthogonality involved in this algorithm?
 }
\end{center}
To answer this question, restrict to the 2-dimensional case and the computation
of the part
\[
\sum_{ t_1' > t_1}  
           a \left(  \sum_{{t_2' \leq t_2}} w_{(t_1',t_2')}, \varphi_{(t_1,t_2),(i_1,i_2)} \right) .
\]
This means that $D=1$ and $d=2$.
Obviously, the whole algorithm would be correct, if we would prolongate all data to a full grid, apply the discretization stencil and then perform restrictions.
Such a prolongation restriction scheme is depicted in Figure \ref{ZZZ_full}.
Yet, this would lead to a computational inefficient algorithm since it requires to perform computations on the complete full grid.
Instead all computations are omitted which are not needed due to the semi-orthogonality property.
To explain this in more detail, consider two overlapping basis functions as in Figure \ref{ZZZ_semi}.
Computations along the algorithmic path depicted in Figure \ref{ZZZ_missing} are needed to take into account the corresponding value in the stiffness matrix $a_n^{\semiortho}(\varphi_{\bt,\bi},\varphi_{\bt',\bi'} )$.
However, by using semi-orthogonality property the result of these computations is zero.
Therefore, this algorithmic path can be omitted. 
This shows that the remaining algorithmic paths depicted in Figure \ref{ZZZ_sparse} are enough for obtaining a correct computation.
This proves that Algorithm \ref{AlgMatrix} correctly calculates a matrix vector multiplication using semi-orthogonality.

\section{The  complete iterative solver}
\vspace{-2pt}

Discretization \ref{DisB} leads to the linear equation system
\begin{equation}
\label{linequSys}
 \cA_n \vec{U} = \cM \vec{F},
\end{equation}
where $\cA_n$ is the stiffness matrix (\ref{StiffnessMat}), $\cM$ the mass matrix, $\vec{F}$ the right hand side in
pre-wavelet format, and $\vec{U}$ the solution vector in pre-wavelet format.
By using the orthogonality property (\ref{orthoPrewave}), $\cM$ is reduced to a diagonal matrix. 

For solving (\ref{linequSys}), the conjugate gradient with a simple diagonal Jacobi precoditioner is applied.
The condition number of $\cA_n$, including this simple preconditioner, is $O(1)$. This follows from the multilevel theory
in \cite{Osw} and the following equivalence of norms:
\[
 \| u \|_{H^1}^2 \cong 
       \sum_{|\bt| \leq n, \bi \in \bXi_\bt } \left| c_{\bt,\bi} \right|^2  a( \varphi_{\bt,\bi},\varphi_{\bt,\bi} ) ,
\]
where $u = \sum_{|\bt| \leq n, \bi \in \bXi_\bt } c_{\bt,\bi}  \varphi_{\bt,\bi}$

\begin{figure}
 \centering
 \includegraphics[width=1.0\textwidth]{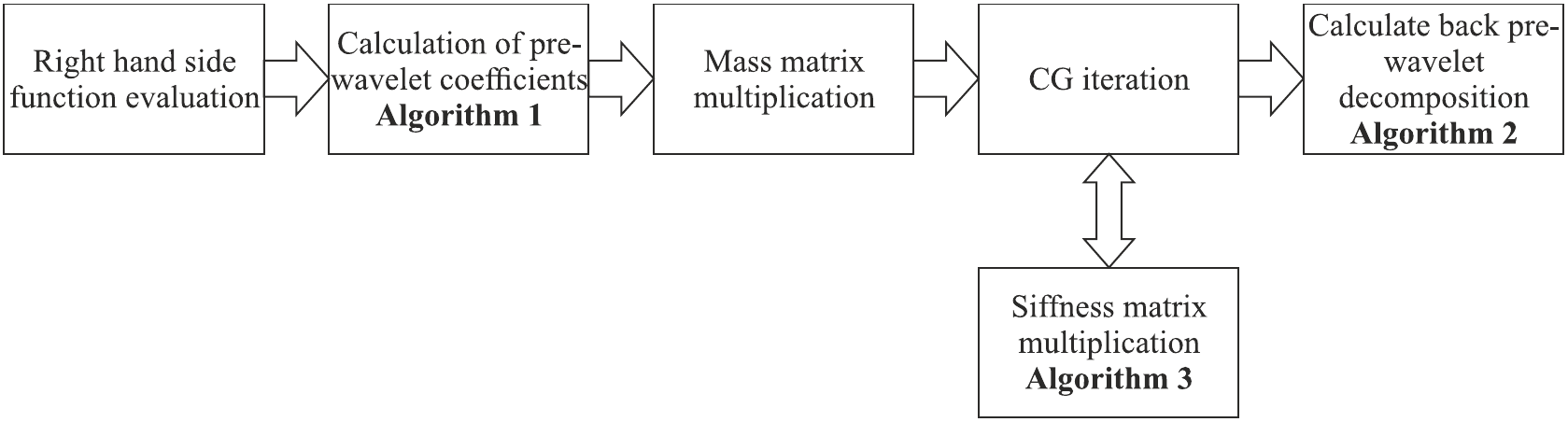}
 \caption{Preconditioned CG-solver using pre-wavelet decomposition.}
 \label{fig_solver}
\end{figure}

Algorithm \ref{prew_coeff} is used to calculate the right hand side vector $\vec{F}$ 
for a  given right hand side vector $(f(p))_{p \in \cD_n}$. The conjugate gradient algorithm requires the application the stiffness matrix multiplication
Algorithm \ref{AlgMatrix}. The resulting vector $\vec{U}$ applied to Algorithm \ref{prew_back}
leads to an approximation of the 
finite element solution $u_{\cD_n}^\sPre \in V_{\cD_n}^\sPre$ of Discretization \ref{DisB}.
The total algorithm is depicted in Figure \ref{fig_solver}.

\section{Numerical results}
\label{simResults}
To show the efficiency of the discretization in this paper, two numerical results are presented.
The first example shows that sparse grids can be used to discretize elliptic partial differential equations
on curvilinear bounded domains in 3D. The second example is a 6-dimensional Helmholtz problem with a variable
coefficient. To our knowledge, this is the first Ritz-Galerkin finite element discretization 
of an elliptic PDE with variable coefficients in a high dimensional space. 

In this paper the discretization stencils were obtained by analytic calculations. If this is not possible,
then one has to interpolate the variable coefficients by a piecewise
constant interpolation of the variable coefficients on the sparse grid as in \cite{dorpfl}. However, this
paper was restricted to analytic calculation of the $27$-stencils and $729$-stencils.

The simulation results were obtained on a workstation without parallel computing. The algorithms
were implemented in C++ programming language.

\begin{figure}[hbt]
 \centering
 \setcounter{subfigure}{0}
 \subfigure[error in $L_{2}$-norm]{
   \label{sparsegrid3D_L2}
   \includegraphics[width=0.48\textwidth]{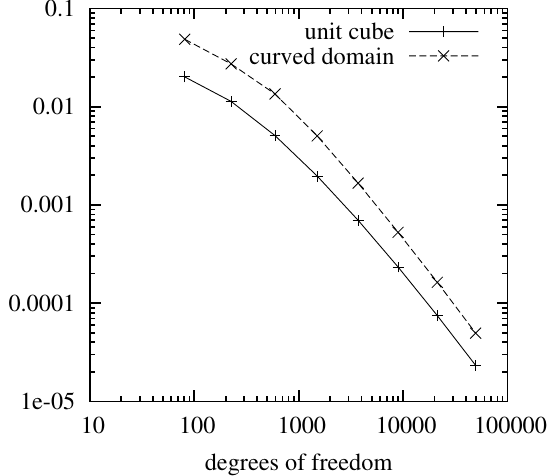}
 }
 \subfigure[error in $L_{\infty}$-norm]{
   \label{sparsegrid3D_Linf}\includegraphics[width=0.48\textwidth]{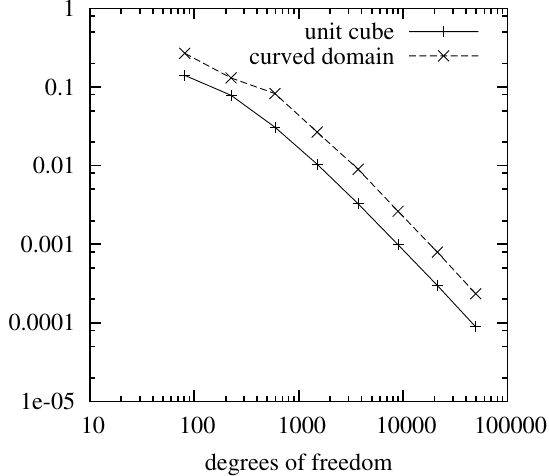}
 }
 \caption{Convergence of the error norm measured by $L_2$ (a) and $L_{\infty}$ (b).}
 \label{sparsegrid3D_error}
\end{figure}


\subsection{Poisson's equation}
\begin{figure}[hbt]
 \centering
 \setcounter{subfigure}{0}
 \subfigure[unit cube]{
   \label{sparsegrid3D}\includegraphics[width=0.45\textwidth]{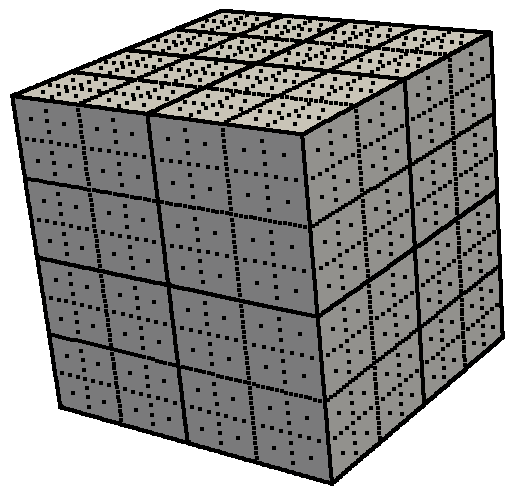}
 }
 \subfigure[curved edges]{
   \label{sparsegrid3D_trans}\includegraphics[width=0.45\textwidth]{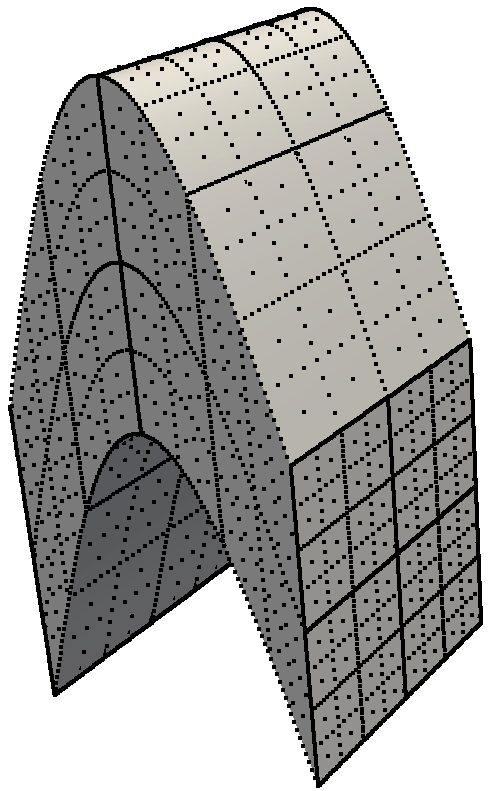}
 }
 \caption{Sparse grid with depth $ n = 7 $ on 3-dimensional unit cube (a) and domain with curved edges (b).}
\end{figure}

We want to show that sparse grids can be applied to discretize partial differential equations
on curvilinear bounded domains. To this end consider Poisson's equation
\begin{eqnarray}
 \label{PoiExE}
 -\Delta u & = & f \; \text{in} \; \Omega \subset \real^3 \\
 u & = & g  \; \text{on} \; \partial \Omega.
 \nonumber
\end{eqnarray}
In order to discretize this problem on a curvilinear bounded domain, one has to subdivide
the domain $\Omega$ into several blocks, such that each of these blocks can smoothly be transformed
to a unit cube. By these transformations to a unit cube, one obtains partial differential equations with variable coefficients on 
the unit cube. The basic idea of this concept is explained in \cite{dorpfl,gordon93} for 2-dimensional domains.
To show that this concept can be extended to a 3-dimensional domain, it is applied to
the curvilinear bounded domain depicted in Figure \ref{sparsegrid3D_trans}. Simulation results
are compared with simulations
obtained on a simple cubical domain $ \left[0,1 \right]^3 $ (see Figure \ref{sparsegrid3D}).
Right hand side $f$ and inhomogeneous Dirichlet boundary condition $g$ are chosen such that 
\[
 u = \sin\left(x \pi \right) \sin\left(y \pi \right) \sin\left(z \pi \right)
\]
is the  exact solution of (\ref{PoiExE}). The curvilinear bounded domain is obtained by the
transformation of the x-coordinate according:
\[
 \tilde{x} = x + \sin\left( y \pi \right) - \frac{1}{2}.
\]
This analytic transformation allows an analytic calculation of the 27-stencils on each subgrid of the sparse grid.

Now, let $u$ be the exact solution of Poisson's problem (\ref{PoiExE})
and $ u_{\cD_n}^\sPre $ the finite element discretization (\ref{equDisSemi})
on the sparse grid $\cD_n$ of depth $n$. Furthermore, let
$e_{n,\infty} = \| u - u_{\cD_n}^\sPre \|_{\infty} $ be the error in the maximum norm
and $ e_{n,2} = \| u - u_{\cD_n}^\sPre \|_{2} $ the error in a suitable weighted discrete $L_2$-norm such 
that a constant error has an equally error norm with respect to the  discrete $L_2$-norm.
Table \ref{tabLInf}, Table \ref{tabL2}, and Figure \ref{sparsegrid3D_error} show that the discretization 
with semi-orthogonality and prewavelets leads to an optimal convergence according the approximation properties of sparse grids.
Moreover, the condition number of the stiffness matrix using a simple diagonal preconditioner stays
below $10$ for $n=2,...,9$ (see Table \ref{tabCond}). Therefore, only a few cg-iterations are needed to
obtain a small algebraic error.

\begin{table}[hbt]
  \centering
  \begin{tabular}{|c | r | r | r | r | r | r | r |}\hline
    & & \multicolumn{2}{c|}{unit cube}& \multicolumn{2}{c|}{curved domain} \\
    $ t_{\text{max}} $ & DOF & \multicolumn{1}{c|}{$ e_{n,\infty} $} & \multicolumn{1}{c|}{$ \frac{e_{n,\infty}}{e_{n-1,\infty}} $} & \multicolumn{1}{c|}{$ e_{n,\infty} $} & \multicolumn{1}{c|}{$ \frac{e_{n,\infty}}{e_{n-1,\infty}} $} \\\hline 
    2 &     7 & $1.4131\e{-1}$ &      & $1.4131\e{-1}$ &      \\
    3 &    31 & $7.8379\e{-2}$ & 1.80 & $1.3066\e{-1}$ & 2.06 \\
    4 &   111 & $3.0813\e{-2}$ & 2.54 & $8.2763\e{-2}$ & 1.58 \\
    5 &   351 & $1.0418\e{-2}$ & 2.96 & $2.6655\e{-2}$ & 3.10 \\
    6 &  1023 & $3.2704\e{-3}$ & 3.19 & $8.9742\e{-3}$ & 2.97 \\
    7 &  2815 & $9.8289\e{-4}$ & 3.33 & $2.6322\e{-3}$ & 3.41 \\
    8 &  7423 & $2.9844\e{-4}$ & 3.29 & $7.9659\e{-4}$ & 3.31 \\
    9 & 18943 & $8.8898\e{-5}$ & 3.36 & $2.3461\e{-4}$ & 3.40 \\\hline
  \end{tabular} 
  \caption{$L_{\infty}$-norm of discretization error.}
  \label{tabLInf}
\end{table}

\begin{table}[hbt]
  \centering
  \begin{tabular}{|c | r | r | r | r | r | r | r |}\hline
    & & \multicolumn{2}{c|}{unit cube}& \multicolumn{2}{c|}{curved domain} \\
    $ t_{\text{max}} $ & DOF & \multicolumn{1}{c|}{$ e_{n,2} $} & \multicolumn{1}{c|}{$ \frac{e_{n,2}}{e_{n-1,2}} $} & \multicolumn{1}{c|}{$ e_{n,2} $} & \multicolumn{1}{c|}{$ \frac{e_{n,2}}{e_{n-1,2}} $} \\\hline 
    2 &     7 & $2.0150\e{-2}$ &        & $4.8651\e{-2}$ &      \\
    3 &    31 & $1.1328\e{-2}$ & 1.78 & $2.7337\e{-2}$ & 1.78 \\
    4 &   111 & $5.0768\e{-3}$ & 2.23 & $1.3524\e{-2}$ & 2.02 \\
    5 &   351 & $1.9677\e{-3}$ & 2.58 & $5.0468\e{-3}$ & 2.68 \\
    6 &  1023 & $6.9732\e{-4}$ & 2.82 & $1.6640\e{-3}$ & 3.03 \\
    7 &  2815 & $2.3303\e{-4}$ & 3.00 & $5.2647\e{-4}$ & 3.16 \\
    8 &  7423 & $7.4714\e{-5}$ & 3.12 & $1.6289\e{-4}$ & 3.23 \\
    9 & 18943 & $2.3233\e{-5}$ & 3.22 & $4.9463\e{-5}$ & 3.29 \\\hline
  \end{tabular} 
  \caption{$L_{2}$-norm of discretization error.}
  \label{tabL2}
\end{table}

\begin{table}[hbt]
  \centering
  \begin{tabular}{|c | r | r | r | r | r | r | r |}\hline
    & & \multicolumn{2}{c|}{unit cube}& \multicolumn{2}{c|}{curved domain} \\
    $ t_{\text{max}} $ & DOF & $ \kappa\left( A \right) $ & $ \kappa\left( C^T A C \right) $ & $ \kappa\left( A \right) $ & $ \kappa\left( C^T A C \right) $   \\\hline 
    2 &     7 &   2.96 & 1.62 &   4.47 & 2.18 \\
    3 &    31 &   8.48 & 2.42 &  23.54 & 3.25 \\
    4 &   111 &  18.82 & 3.87 &  64.25 & 3.91 \\
    5 &   351 &  48.57 & 4.40 & 121.38 & 4.49 \\
    6 &  1023 & 126.90 & 4.59 & 179.78 & 4.99 \\
    7 &  2815 & 207.40 & 5.05 & 232.33 & 5.06 \\
    8 &  7423 & 267.92 & 5.10 & 270.43 & 5.43 \\
    9 & 18943 & 283.70 & 5.46 & 306.03 & 5.95 \\\hline
  \end{tabular} 
  \caption{Condition number of pre-wavelet discretization and diagonal preconditioning.}
  \label{tabCond}
\end{table}

\subsection{Helmholtz equation with variable coefficients in high dimensions}
Consider the 6-dimensional Helmholtz problem
\begin{eqnarray}
\label{eqnHelmsix}
 -\Delta u + c u & = & f \; \text{in} \; \Omega := \left[ 0,1 \right]^6 \\
 \nonumber
 u & = & 0  \; \text{on} \; \partial \Omega
\end{eqnarray}
with variable coefficient
\begin{eqnarray}
\nonumber
\lefteqn{
 c\left( x,y,z,u,v,w \right) :=} \\
 \label{varCoeff}
&& \left(1-x^2\right)\left(1-y^2\right)\left(1-z^2\right)\left(1-u^2\right)\left(1-v^2\right)\left(1-w^2\right) .
\end{eqnarray}
The right hand side  $f$ is chosen such that 
\[
u =  \sin\left(x \pi \right) \sin\left(y \pi \right) \sin\left(z \pi \right) \sin\left(u \pi \right) \sin\left(v \pi \right) \sin\left(w \pi \right)
\]
is a solution of (\ref{eqnHelmsix}).

Table \ref{tab6d} shows the discretization error for constant coefficients ($c = 1$) as well as 
for the variable coefficient (\ref{varCoeff}).
The convergence rate for the problem with constant coefficients is similar to the reported
convergence behavior in \cite{ZengerSparseGrids}. 
In addition, the solution of the discretization with semi-orthogonality and prewavelets 
convergences in case of variable coefficients is as fast as in case of constant coefficients.
This shows that the discretization with semi-orthogonality does not introduce any remarkable additional 
errors. 

\begin{table}[hbt]
  \centering
  \begin{tabular}{|c | r | r | r | r | r |}\hline
     &  & \multicolumn{2}{c|}{constant coefficient} & \multicolumn{2}{c|}{variable coefficient} \\
    $ t_{\text{max}} $ & DOF & $ e_{n,\infty} $ & \multicolumn{1}{c|}{$ \frac{e_{n,\infty}}{e_{n-1,\infty}} $} & $ e_{n,\infty} $ & \multicolumn{1}{c|}{$ \frac{e_{n,\infty}}{e_{n-1,\infty}} $} \\ \hline 
    2 &    13 & $ 0.40427 $ &      & $ 0.391625 $ &      \\
    3 &    97 & $ 0.28396 $ & 1.42 & $ 0.282802 $ & 1.38 \\
    4 &   545 & $ 0.10686 $ & 2.65 & $ 0.106861 $ & 2.64 \\
    5 &  2561 & $ 0.04028 $ & 2.65 & $ 0.040046 $ & 2.66 \\
    6 & 10625 & $ 0.01533 $ & 2.62 & $ 0.015328 $ & 2.61 \\\hline
  \end{tabular} 
  \caption{Error convergence for a 6-dimensional Helmholtz problem.}
  \label{tab6d}
\end{table}

\section{Acknowledgment}
\vspace{-2pt}
The authors gratefully acknowledge funding of the Erlangen Graduate School in Advanced Optical Technologies (SAOT)
by the German Research Foundation (DFG) in the framework of the German excellence initiative.

\bibliographystyle{wileyj}
\bibliography{lithere}

\end{document}